\documentclass[runningheads]{llncs}
\usepackage{amsfonts}
\usepackage{epsfig}
\usepackage{amsmath}
\usepackage{graphicx}
\usepackage{amsmath}
\usepackage{amscd}
\usepackage{stmaryrd}
\usepackage{algorithm}
\usepackage{arydshln}
\usepackage{algorithmic}
\usepackage{youngtab}
\usepackage{young}
\usepackage{geometry}
\usepackage{bbm}
\usepackage{leftidx}

\usepackage[colorlinks,
            linkcolor=blue,
            anchorcolor=blue,
            citecolor=blue,
            ]{hyperref}

\def\bt{\begin{theorem}}
\def\et{\end{theorem}}
\def\bp{\begin{proposition}}
\def\ep{\end{proposition}}
\def\bc{\begin{corollary}}
\def\ec{\end{corollary}}
\def\bo{\begin{proof}}
\def\eo{\end{proof}}
\def\bx{\begin{example}}
\def\ex{\end{example}}
\def\br{\begin{remark}}
\def\er{\end{remark}}
\def\bl{\begin{lemma}}
\def\el{\end{lemma}}

\def\bn{\begin{definition}}
\def\en{\end{definition}}
\def\ba{\begin{array}}
\def\ea{\end{array}}
\def\be{\begin{equation}}
\def\ee{\end{equation}}
\def\bd{\begin{description}}
\def\ed{\end{description}}
\def\bu{\begin{enumerate}}
\def\eu{\end{enumerate}}
\def\bi{\begin{itemize}}
\def\ei{\end{itemize}}

\newbox\bigstrutbox
\setbox\bigstrutbox=\hbox{\vrule height12.5pt depth5pt width0pt}
\def\bigstrut{\relax\ifmmode\copy\bigstrutbox\else\unhcopy\bigstrutbox\fi}

\newbox\Bigstrutbox
\setbox\Bigstrutbox=\hbox{\vrule height17.5pt depth5pt width0pt}
\def\Bigstrut{\relax\ifmmode\copy\Bigstrutbox\else\unhcopy\Bigstrutbox\fi}

\def\ds{\displaystyle}

\def\A{{\bf A}}
\def\B{{\bf B}}
\def\C{{\bf C}}
\def\D{{\bf D}}

\def\R{{\bf R}}

\def\a{{\bf a}}
\def\b{{\bf b}}
\def\c{{\bf c}}
\def\d{{\bf d}}
\def\e{{\bf e}}
\def\f{{\bf f}}
\def\g{{\bf g}}

\def\m{{\bf m}}

\def\r{{\bf r}}
\def\s{{\bf s}}
\def\t{{\bf t}}
\def\u{{\bf u}}
\def\v{{\bf v}}
\def\w{{\bf w}}
\def\x{{\bf x}}
\def\y{{\bf y}}

\def\0{{\bf 0}}
\def\1{{\bf 1}}
\def\2{{\bf 2}}
\def\3{{\bf 3}}
\def\4{{\bf 4}}
\def\5{{\bf 5}}
\def\6{{\bf 6}}
\def\7{{\bf 7}}
\def\8{{\bf 8}}
\def\9{{\bf 9}}

\begin{document}

\pagestyle{headings}

\mainmatter

\title{An Improvement of Rota's Straightening Algorithm}

\titlerunning{An Improvement of Rota's Straightening Algorithm}

\author{Changpeng Shao \\
cpshao@amss.ac.cn}

\authorrunning{Changpeng Shao}

\institute{Key Laboratory of System Control, Academy of Mathematics and Systems Science, \\
Chinese Academy of Sciences, Beijing 100190, China}

\maketitle

\begin{abstract}
In bracket algebra, the calculation of invariant division and invariant Gr\"{o}bner basis proposed in \cite{li 2014} rely on straightening algorithm. Until now, there are at least three different types of straightening algorithms, among which Rota's straightening algorithm has the best efficiency. However, there exists a flaw in Rota's straightening algorithm, i.e., it needs find all the straight bracket monomials with the same content as the input beforehand, which is quite difficult. So in this paper, we will propose a new straightening algorithm based on dual bracket, which is a new concept of Young tableau. In this new straightening algorithm, we only need to find a few number of straight bracket monomials in each step instead of finding them all in one step. And so it is an improvement of Rota's straightening algorithm. According to our tests, this new straightening algorithm reflects more advantages when the dimension and the degree increase. Moreover, this straightening algorithm still works when Rota's straightening algorithm fails.
\end{abstract}


\section{Introduction}
\setcounter{equation}{0}

The main study objects of classical invariant theory are invariants of the special linear group. Given a $n$ dimensional $\mathbb{K}$-vector space $\mathcal{V}$, and $m$ vector variables
\be
\a_i=(x_{i1},x_{i2},\cdots,x_{in})^T, \hspace{.4cm} i=1,2,\cdots,m
\ee
of $\mathcal{V}$. The first fundamental theorem of classical invariant theory \cite{li}, \cite{sturmfels 2} states that all $SL(n)$-invariants lie in the quotient ring of the polynomial ring generated by determinants
\be \label{all determinants}
\{[\a_{i_1}\a_{i_2}\cdots\a_{i_n}]:=\det(x_{i_rs})_{n\times n}\mid1\leq i_1<i_2<\cdots<i_n\leq m\}.
\ee
This polynomial quotient ring is called \emph{bracket algebra}, and determinant is called \emph{basic invariant}. Polynomials in bracket algebra are called bracket polynomials.

In bracket algebra, each determinant is a variable. However, if expanding all these determinants and considered as polynomials in another polynomial ring
\be
\mathbb{K}[\{x_{ij}\mid i=1,2,\cdots,m;j=1,2,\cdots,n\}],
\ee
we will find that: There exist syzygies among determinants, i.e., there exist non bracket polynomials which equal to zero after expanding all determinants. The ideal generated by all syzygies in the polynomial ring generated by determinants (\ref{all determinants}) will be denoted as $\mathcal{I}$, and called the \emph{syzygy ideal} of bracket algebra.

A bracket polynomial can be written in its normal form via Young's straightening algorithm \cite{rutherford}, \cite{young03}, \cite{young}, this procedure is called \emph{straightening}. From the viewpoint of computer algebra \cite{sturmfels 2}, Young's straightening algorithm is equivalent to the reduction procedure of a Gr\"{o}bner bases of $\mathcal{I}$ \cite{sturmfels and white}.

The invariant division in bracket algebra proposed in \cite{li 2014} is an algebraic operation among straight bracket polynomials, so we need straighten each bracket polynomial in the process of this invariant division. Classical straightening algorithm \cite{young03}, \cite{young} was proposed by Young in 1928. Later, in 1991, White \cite{white1} examined the choices that are involved in implementing a straightening algorithm and described some particular variations which they found to be relatively efficient, then gave five different straightening algorithms. These five straightening algorithms all possess a high efficiency than Young's straightening algorithm based on his tests and our tests in section 5. These six straightening algorithms are reduction procedures of different Gr\"{o}bner bases of the syzygy ideal $\mathcal{I}$.

In 1974, Doubilet, Rota, and Stein \cite{rota-foundation} initially studied the theory structure about invariant theory in the case of character zero and provided a straightening algorithm. Later, in 1978, D\'{e}sarm\'{e}nien, Kung, and Rota \cite{rota-straightening-formula} considered this problem furthermore, and using Capelli formula \cite{capelli} gave a self-contained combinatorial presentation about the validity of the straightening algorithm in \cite{rota-foundation}. Rota's straightening algorithm was also used to give a simple proof of the first and second fundamental theorem of classical invariant theory. Different from the above two types of straightening algorithms, Rota's straightening algorithm is a procedure of the solving of a system of linear equations.

By comparison, Rota's straightening algorithm can deal with bracket polynomials of any degree in each step. However, White's straightening algorithms and the classical straightening algorithm can only deal with bracket polynomial of degree two in each step. The most common case, however, in practice has degree larger than two. So from this point, Rota's straightening algorithm depicts more advantages.

On the other hand, there is an awe-inspiring step in Rota's straightening algorithm, i.e., find all straight bracket monomials with the same content as the input bracket polynomial. This certainly can only be achieved by computer. Maybe because of this, few people prefer to use Rota's straightening algorithm to straighten bracket polynomials. From a lot of tests, we find that:
\begin{enumerate}
  \item When the degree and dimension are relatively small, the enumeration of all straight bracket monomials is not very difficult.
  \item Assume that all straight bracket monomials have been got, then the texting results show that Rota's straightening algorithm is much better than White's straightening algorithm. The reason in my view is: Rota's straightening algorithm is a procedure of the solving of a system of linear equations $\A\x=\b$, where the coefficient matrix $\A$ is a lower triangular matrix with diagonal elements nonzero. Matrix $\A$ depends on the degree, dimension and the content of the input bracket polynomial. So, in the straightening process, if the degree, dimension and the content keep invariant, then $\A$ is a common coefficient matrix, and only needed computed once.
  \item When the degree and dimension are larger than five, the enumeration of all straight bracket monomials becomes difficult. If the enumeration of all straight bracket monomials fails, then Rota's straightening algorithm can not works.
\end{enumerate}

Therefore, in order to improve the efficiency of straightening algorithm and invariant division, we need to make improvements on Rota's straightening algorithm.

In \cite{li 2014}, we established the relationship between bracket algebra and its coordinate polynomial ring. One result shows that there exists an unique one to one correspondence between straight monic bracket monomials and standard monic coordinate monomials. From this correspondence, it's not too hard to give a straightening algorithm in bracket algebra from the viewpoint of coordinate polynomial ring.

From the relationship between bracket algebra and its coordinate polynomial ring, we introduced the concept of dual bracket. Dual bracket can be written in the Young tableau form just like bracket monomials. Given a Young tableau, compress each column by Hadamard product as a vector, then using permanent or determinant of matrices to define two types of dual bracket. A main result about dual bracket is that: Each bracket monomial is a polynomial of dual brackets. Thus, this gives us a method to study bracket monomials from the viewpoint of dual bracket.

Since dual bracket and coordinate monomial also has an one to one correspondence just like the one to one correspondence between straight monic bracket monomials and standard monic coordinate monomials, so the idea of giving a straightening algorithm in bracket algebra based on coordinate polynomial ring is easily extended to dual bracket. However, these two straightening algorithms do not have a high efficiency. The main reason relies on that in the expanding formula of bracket monomial into dual bracket, not every term is needed. Only some special dual bracket is necessary. Combining this special terms together, we will get a straight bracket polynomial, denoted as $SB(f)$, where $f$ is a bracket polynomial. Therefore, according to $SB(f)$, we get the new straightening algorithm in bracket algebra.

The essence of this new straightening algorithm and the straightening algorithm based on coordinate polynomial ring are the same. The main advantage of this new straightening algorithm is: Dual bracket contains the same form (i.e., Young tableau form) as bracket monomials, so in the calculation of $SB(f)$, we may not need to expand bracket monomial into dual brackets. This is impossible for coordinate polynomial ring. Note that the expanding expression contains $(n!)^{d-1}$ terms, where $n$ is the dimension, and $d$ is the degree of the bracket monomial. So, we can see that expanding will decrease the efficiency of the straightening algorithm.

In this paper, we will also present an algorithm of computing $SB(f)$, which is independent of its definition. The elementary idea of this algorithm comes from the algorithm of finding all straight bracket monomials. In the algorithm of computing $SB(f)$, there are four implied rules. Deleting the third rule is the algorithm of finding all straight bracket monomials. This figures out that, we use the third rule to give a restriction of finding all straight bracket monomials in Rota's straightening algorithm. Thus we only need to compute a few straight bracket monomials in each step. At the end of this new straightening algorithm, the union of all terms of $SB(f)$ is approximate to all straight bracket monomials. Hence, the difficulty in Rota's straightening algorithm is actually solved by parts in the new straightening algorithm based on $SB(f)$.

Form a lot of tests, we find that our straightening algorithm is much better than White's straightening algorithms, while a little worse than Rota's straightening algorithm if all the straight bracket monomials can be enumerated. However, when the degree of bracket monomials are higher than $5$, all straight bracket monomials are not easy to get anymore, and thus Rota's straightening algorithm can not works any more. But our straightening algorithm still works, so from this point, this is a big improvement of Rota's straightening algorithm. On the other hand, the results about straightening coefficients obtained in \cite{clausen2}, \cite{desarmenien}, \cite{huang} based on Rota's straightening algorithm can also be achieved by our straightening algorithm with a much more simple proof (Appendix II).

This paper is organized as follows. In section 2, we will review some basic definitions about Young tableau and bracket algebra and some necessary results from the paper \cite{li 2014}. In section 3, we will briefly introduce the classical straightening algorithm, White's straightening algorithm and Rota's straightening algorithm. Section 4 is the main parts of this paper, it is devoted to give the definition of dual bracket, present some properties of dual bracket and the new straightening algorithm. In section 5, we will make an analysis on the efficiency of the the three famous straightening algorithm and our straightening algorithm via specific examples.

\textbf{Notation}: in this paper, all the bold letters such as $\b_i,\c_i,\d_i,\y_i,\a_{ij},\b_{ij},\c_{ij},\d_{ij},\f_{ij},\g_{ij},\y_{ij},\ldots$ except $\e_i,\r_i$ without other description refers to elements of $\{\a_1,\a_2,\ldots,\a_m\}$.

\section{Young tableau and bracket algebra}
\setcounter{equation}{0}

In this section, we will give a brief introduction about Young tableau and bracket algebra \cite{rota-straightening-formula}, \cite{fulton}, \cite{li}, and will introduce some concepts about Young tableau for the convenience of the following sections.

Let $\mathcal{V}$ be a $n$ dimensional $\mathbb{K}$-vector space, where char$(\mathbb{K})\neq2$. Choosing $m\geq n$ vector variables $\a_1,\a_2,\cdots,\a_m$ from $\mathcal{V}$, where the coordinate of $\a_i (1\leq i\leq m)$ is
\be
\a_i=(x_{i1},x_{i2},\cdots,x_{in})^T.
\ee

\bn \cite{li}
The bracket of $\a_{i_1},\a_{i_2},\cdots,\a_{i_n}$ is denoted as $[\a_{i_1}\a_{i_2}\cdots\a_{i_n}]$ and defined as
\be
[\a_{i_1}\a_{i_2}\cdots\a_{i_n}]:=\left|
                            \begin{array}{cccc}
                              x_{i_11} & x_{i_21} & \cdots & x_{i_n1} \\
                              x_{i_12} & x_{i_22} & \cdots & x_{i_n2} \\
                              \vdots   & \vdots   & \ddots & \vdots \\
                              x_{i_1n} & x_{i_2n} & \cdots & x_{i_nn} \\
                            \end{array}
                          \right|.
\ee
The $n$ dimensional bracket algebra generated by $\a_1,\a_2,\cdots,\a_m$, denoted as $\mathcal{B}[\{\a_i:i=1,2,\cdots,m\}]$, is the $\mathbb{K}$-polynomial ring:
\be
\mathbb{K}[\{[\a_{i_1}\a_{i_2}\cdots\a_{i_n}]\mid1\leq i_1,i_2,\cdots,i_n\leq m\}].
\ee
module the ideal, which is called the syzygy ideal of bracket algebra, generated by:
\begin{description}
  \item[\hspace{.2cm} B1.] \hspace{.2cm} $[\a_{i_1}\a_{i_2}\cdots\a_{i_n}]$ if $i_j=i_k$ for some $j\neq k$.
  \item[\hspace{.2cm} B2.] \hspace{.2cm} $[\a_{i_1}\a_{i_2}\cdots\a_{i_n}]-\emph{sign}(\sigma)[\a_{i_{\sigma(1)}}\a_{i_{\sigma(2)}}\cdots\a_{i_{\sigma(n)}}]$ for any permutation $\sigma$ of $1,2,\cdots,n$.
  \item[\hspace{.2cm} GP.] \hspace{.2cm} $\ds\sum_{k=1}^{n+1}(-1)^k [\a_{i_1}\a_{i_2}\cdots\a_{i_{k-1}}\a_{i_{k+1}}\cdots\a_{i_{n+1}}][\a_{i_k}\a_{j_1}\a_{j_2}\cdots\cdots\a_{j_{n-1}}]$.
\end{description}
\en

The fundamental notion in bracket algebra is Young tableau \cite{rota-straightening-formula}, \cite{fulton}.
Let $(\lambda)=(\lambda_1,\lambda_2,\cdots,\lambda_s)$ be a partition of the integer $n$: that is, $(\lambda)$ is a finite sequence of positive integers such that
\be
\lambda_1+\lambda_2+\cdots+\lambda_s=n, \hspace{2cm} \lambda_1\geq\lambda_2\geq\cdots\geq\lambda_s>0.
\ee
If $(\lambda)$ is a partition of $n$, its shape, also denoted by $(\lambda)$, is the set of integer points $(i,j)$ in the plane, with $1\leq j\leq s$ and $1\leq i\leq \lambda_j$.

A Young tableau on the shape $(\lambda)$ with values in the set $\mathcal{A}$ is an assignment of an element of $\mathcal{A}$ to each point in the shape $(\lambda)$. For example, Young tableaux of shape $(\lambda)=(5,4,2,2,1,1)$ with values in the integers can be:
\[T_1=\begin{array}{c}
          31265 \\
          5213\hfill \\
          63\hfill \\
          45\hfill \\
          9\hfill \\
          9\hfill \\
        \end{array} \hspace{1cm} T_2=\begin{array}{c}
          22222 \\
          2334\hfill \\
          54\hfill \\
          66\hfill \\
          8\hfill \\
          7\hfill \\
        \end{array}\hspace{1cm} T_3=\begin{array}{c}
          12345 \\
          2345\hfill \\
          34\hfill \\
          45\hfill \\
          6\hfill \\
          7\hfill \\
        \end{array}\hspace{1cm} T_4=\begin{array}{c}
          17899 \\
          2789\hfill \\
          39\hfill \\
          49\hfill \\
          5\hfill \\
          6\hfill \\
        \end{array}\hspace{1cm} \cdots\cdots\]

In this paper, $\mathcal{A}$ is always assumed as $\{\a_1,\a_2,\cdots,\a_m\}$ with order $\a_1\prec \a_2\prec \cdots \prec \a_m$. We only encounter the Young tableau with the following form:
\be \label{Young-tableau}
\begin{array}{ccccc}
        \y_{11} & \y_{12} & \cdots  & \y_{1n}   \\
        \y_{21} & \y_{22} & \cdots  & \y_{2n}   \\
        \vdots  & \vdots  & \ddots  & \vdots    \\
        \y_{k1} & \y_{k2} & \cdots  &  \y_{kn}  \\
      \end{array}
\ee
where $n$ is called the \emph{dimension} of the tableau, $k$ is called the \emph{degree} of the tableau. The following are some necessary definitions about Young tableau in this paper.

\bn Let $Y$ be a Young tableau with the form (\ref{Young-tableau}).

(1). Suppose that $1\leq i\leq k$ and $1\leq j\leq n$, define the $i$-th row of $Y$ as the following Young tableau with degree 1, and dimension $n$:
\be
\y_{i1}~~\y_{i2}~~\cdots~~\y_{in},
\ee
Define the $j$-th column of $Y$ as the following Young tableau with degree $k$, and dimension $1$:
\be
\begin{array}{c}
        \y_{1j}  \\
        \y_{2j}  \\
        \vdots   \\
        \y_{kj}  \\
      \end{array}.
\ee

(2). $Y$ is called straight if the entries in each row are increasing from left to right, and the entries in each column are nondecreasing downward.
\en

\bn Assume that $U$ is a Young tableau of degree $u$, dimension $k$, $V$ is a Young tableau of degree $v$, dimension $l$,

(1). If $u=v$, define the row joint operation $\circ_R$ among them as:
\be
X\circ_RY:=X~~Y.
\ee

(2). If $k=l$, define the column joint operation $\circ_C$ among them as:
\be
X\circ_CY:=\begin{array}{c}
        X \\
        Y \\
      \end{array}.
\ee
\en

\bn Let $X$ be a Young tableau of degree $d$, dimension $n$, define

(1). the row normalization $\mathfrak{N}_R(X)$ of $X$ as: Sorting each row of $X$, and do some row permutations on $X$ such that for any $1\leq i\leq d-1$, the $i$-th row has a lower lexicographical order than the $(i+1)$-th row of $X$ when compared from left to right.

(2). the column normalization $\mathfrak{N}_C(X)$ of $X$ as: Sorting each column of $X$, and do some column permutations on $X$ such that for any $1\leq j\leq n-1$, the $j$-th column has a lower lexicographical order than the $(j+1)$-th column of $X$ when compared from top to bottom.
\en

\bx Let $X=\begin{array}{ccc}
               \a_2 & \a_7 & \a_4 \\
               \a_3 & \a_8 & \a_9 \\
               \a_1 & \a_6 & \a_5 \\
             \end{array}$, then
$\mathfrak{N}_R(X)=\begin{array}{ccc}
               \a_1 & \a_5 & \a_6 \\
               \a_2 & \a_4 & \a_7 \\
               \a_3 & \a_8 & \a_9 \\
             \end{array}$ and
$\mathfrak{N}_C(X)=\begin{array}{ccc}
               \a_1 & \a_4 & \a_6 \\
               \a_2 & \a_5 & \a_7 \\
               \a_3 & \a_9 & \a_8 \\
             \end{array}.$
\ex

\bn \cite{li}
Let
\be \label{bracket monomial}
f=[\a_{11}\a_{12}\cdots\a_{1n}][\a_{21}\a_{22}\cdots\a_{2n}]\cdots[\a_{d1}\a_{d2}\cdots\a_{dn}],
\ee
be a monic monomial in bracket algebra, called bracket monomial, where $\a_{ij}\in \{\a_1,\cdots,\a_m\}$. Setting
\be \label{matrix form}
F=\begin{array}{cccc}
    \a_{11} & \a_{12} & \cdots & \a_{1n} \\
    \a_{21} & \a_{22} & \cdots & \a_{2n} \\
    \vdots  & \vdots  & \ddots & \vdots \\
    \a_{d1} & \a_{d2} & \cdots & \a_{dn} \\
  \end{array}
\ee
which is called the Young tableau of $f$. In addition, define
\be
[F]=\left[\begin{array}{cccc}
    \a_{11} & \a_{12} & \cdots & \a_{1n} \\
    \a_{21} & \a_{22} & \cdots & \a_{2n} \\
    \vdots  & \vdots  & \ddots & \vdots \\
    \a_{d1} & \a_{d2} & \cdots & \a_{dn} \\
  \end{array}\right]=f.
\ee
\en

Due to the anti-symmetry property of determinant, a bracket monomial $f$ in the form (\ref{bracket monomial}) satisfies:
\be\ba{lll} \vspace{.1cm}
&& \a_{i_1}\prec \a_{i_2} \prec \cdots \prec \a_{i_n}, ~~~~ i=1,2,\cdots, d; \\
&& (\a_{i1},\a_{i2},\cdots,\a_{i_n}) \preceq (\a_{(i+1)1},\a_{(i+1)2},\cdots,\a_{(i+1)n}) ~\textmd{lexicographically},~~~~ i=1,2,\cdots, d-1,
\ea\ee
will be called \emph{normal bracket monomial}, or \emph{in normal form}.

\bn  \cite{li} A bracket monomial is said straight if its Young tableau is straight.
\en

\bl  \cite{li}, \cite{weyl} (First Main Theorem of Classical Invariant Theory) Straight bracket monomials comprise a basis of bracket algebra as a $\mathbb{Z}$-module.
\el

The way of writing a bracket polynomial into its normal form is called straightening. Young \cite{young03}, \cite{young} was the first one who proposed an algorithm of straightening. We will return this topic in the next section.

\bn \cite{li 2014}
The negative column order in bracket algebra is defined as: Let $[F], [G]$ be two straight bracket monomials, with degree $d$ and $l$ respectively, if
\begin{enumerate}
  \item $d<l$, or
  \item $d=l$, while there exist $i,j$, such that $F$ and $G$ contain the same first $(j-1)$ column. Denoting the $j$-th column of $F$ and $G$ are
  \[\begin{array}{c}
       \f_{1j} \\
       \f_{2j} \\
       \vdots \\
       \f_{dj} \\
     \end{array}~\textmd{and}~\begin{array}{c}
       \g_{1j} \\
       \g_{2j} \\
       \vdots \\
       \g_{lj} \\
     \end{array},\]
  then for any $k<i$, $\f_{kj}=\g_{kj}$, but $\f_{ij}\succ \g_{ij}$.
\end{enumerate}
Then we call $[F]$ has a low negative column order than $[G]$, denoted as $[F]\prec [G]$.
\en

\bl \cite{li 2014} \label{leading term:bracket monomial}
Let $f$ be a bracket monomial in normal from where each bracket is nonzero. After straightening $f$, the leading term of the normal form in the negative column order is $f^{\downarrow\downarrow}$.
\el

\bn  \cite{li 2014}
The polynomial ring $\mathbb{K}[\{x_{ij}\mid 1\leq i\leq m,1\leq j\leq n\}]$ is called the coordinate polynomial ring of bracket algebra.
\en

There is a natural correspondence between bracket algebra and its coordinate polynomial ring. Any bracket polynomial $f$ after expanding all determinants is a polynomial of the coordinate polynomial ring, denoted as $f^c$. Let $\e_1,\e_2,\cdots,\e_n$ is a basis of $\mathcal{V}$ which satisfy $[\e_1\e_2\cdots\e_n]=1$. Choosing the order:
\[\e_1\prec \e_2\prec \cdots\prec \e_n\prec\a_1\prec \a_2\prec \cdots\prec \a_m.\]
\bn \cite{li 2014}
The degree-lex order generated by the following order, reading from left to right and from top to bottom, in the coordinate polynomial ring is called negative basis order:
\be \label{negative basis order}
\begin{tabular}{lllll}
            & $x_{1n}$      ~~& $\succ x_{2n}$      ~~& $\succ\cdots$   ~~& $\succ x_{mn}$      \\
    $\succ$ & $x_{1(n-1)}$  ~~& $\succ x_{2(n-1)}$  ~~& $\succ\cdots$   ~~& $\succ x_{m(n-1)}$  \\
            & $\vdots$      ~~& $\vdots$            ~~& $\vdots$        ~~& $\vdots$            \\
    $\succ$ & $x_{11}$      ~~& $\succ x_{21}$      ~~& $\succ\cdots$   ~~& $\succ x_{m1}$      \\
\end{tabular}
\ee
\en

\bn
Let $f$ be a bracket monomial, define
\be\ba{lll} \vspace{.1cm}
  LT(f) &=& \textmd{the~leading~term~of}~f~\textmd{under~the~negative~column~order}, \\ \vspace{.1cm}
  LM(f) &=& \textmd{the~leading~monomial~of}~f~\textmd{under~the~negative~column~order}, \\ \vspace{.1cm}
  LC(f) &=& \textmd{the~leading~coefficient~of}~f~\textmd{under~the~negative~column~order}, \\ \vspace{.1cm}
  lt(f^c) &=& \textmd{the~leading~term~of}~f^c~\textmd{under~the~negative~basis~order}, \\ \vspace{.1cm}
  lm(f^c) &=& \textmd{the~leading~monomial~of}~f^c~\textmd{under~the~negative~basis~order}, \\
  lc(f^c) &=& \textmd{the~leading~coefficient~of}~f^c~\textmd{under~the~negative~basis~order}.
\ea\ee
\en
\bn \cite{li 2014}
A monomial $M$ has the following form
\be\ba{lll} \label{tableau-form-coordinate-monomial}
&& x_{i_1^11}~x_{i_1^21}~\ldots ~x_{i_1^d1} \\
&& x_{i_2^12}~x_{i_2^22}~\ldots ~x_{i_2^d2} \\
&&  ~~~~~~~~~\vdots    \\
&& x_{i_n^1n}~x_{i_n^2n}~\ldots ~x_{i_n^dn}
\ea\ee
is called standard monic coordinate monomial, where
\be
1\leq i_1^j<i_2^j<\ldots<i_n^j\leq n,~~~j=1,2,\cdots,d.
\ee
\en
\bl \cite{li 2014} \label{relationship between bracket algebra and its coordinate ring}
Any standard monic coordinate monomial $M$ has the form (\ref{tableau-form-coordinate-monomial}) corresponds to an unique straight Young tableau
\be \label{Young tableau of M}
F=\begin{array}{cccc}
        \a_{i_1^1} & \a_{i_2^1} & \cdots & \a_{i_n^1} \\
        \a_{i_1^2} & \a_{i_2^2} & \cdots & \a_{i_n^2} \\
        \vdots     & \vdots     & \ddots & \vdots     \\
        \a_{i_1^d} & \a_{i_2^d} & \cdots & \a_{i_n^d} \\
      \end{array},
\ee
such that $lt([F]^c)$ is (\ref{tableau-form-coordinate-monomial}).
\el
\bn
The Young tableau (\ref{Young tableau of M}) is called the Young tableau of $M$, denoted as $Tableau(M)$.
\en

\section{Brief introduction of straightening algorithms in bracket algebra}
\setcounter{equation}{0}

In this section, we will briefly introduce three straightening algorithms: the classical straightening algorithm \cite{li}, White's straightening algorithm \cite{white1} and especially Rota's straightening algorithm \cite{rota-straightening-formula}, \cite{rota-foundation}.

\subsection{Classical straightening algorithm}
The classical straightening algorithm in bracket algebra was proposed by Young \cite{young03}, \cite{young}. Let
\be\ba{lll} \vspace{.1cm}
\B_r       &=& \b_1\b_2\cdots\b_r; \\ \vspace{.1cm}
\C_{n+1}   &=& \c_1\c_2\cdots\c_{n+1}; \\
\D_{n-1-r} &=& \d_1\d_2\cdots\d_{n-1-r}
\ea\ee
be three sequences of $\a$'s with length $r,n+1$ and $n-1-r$ respectively and
\be\ba{lll} \vspace{.1cm}
&& \b_1\prec\b_2\prec\cdots\prec\b_r; \\ \vspace{.1cm}
&& \c_1\prec\c_2\prec\cdots\prec\c_{n+1}; \\
&& \d_1\prec\d_2\prec\cdots\prec\d_{n-1-r}.
\ea\ee
The \emph{van der Waerden (VW) syzygy} is
\be \label{vw}
\sum_{(n-r,r+1)\vdash \C_{n+1}}[\B_r\C_{n+1(1)}][\C_{n+1(2)}\D_{n-1-r}]=0,
\ee
where $\C_{n+1(1)},\C_{n+1(2)}$ are the Sweedler notation of two subsequences of length $r+1,n-r$ obtained by bipartitioning $\C_{n+1}$, with the sign of permutation of the partition included in the notation.

Let $f=\left[
       \begin{array}{lllllll}
         \a_{i_1} & \cdots & \a_{i_r} & \a_{j_{r+2}} & \a_{j_{r+3}} & \cdots & \a_{j_{n+1}} \\
         \a_{j_1} & \cdots & \a_{j_r} & \a_{j_{r+1}} & \a_{k_1}     & \cdots & \a_{k_{n-r-1}} \\
       \end{array}
     \right]$ be a non-straight bracket monomial of degree 2, with
\[\a_{i_u} \preceq \a_{j_u},~u=1,\cdots,r; \hspace{1cm} \a_{j_{r+1}} \prec \a_{j_{r+2}}.\]
Then in (\ref{vw}), choosing
\[\B_r=\a_{i_1}\cdots\a_{i_r}, \hspace{1cm} \C_{n+1}=\a_{j_1}\cdots\a_{j_{n+1}}, \hspace{1cm} \D_{n-r-1}=\a_{k_1}\cdots\a_{k_{n-r-1}}.\]
We will get an equivalent expression of $f$ with the first $(r+1)$ column straight. Continue this procedure, and finally we can straighten $f$.

\subsection{White's implementation of classical straightening algorithm}

Besides the VW relations, White also used the \emph{multiple sygyzy}: Let $\C=\c_1\cdots\c_n$ be a sequence about $\a$'s with length $n$, $\{\b_1,\cdots,\b_n\}\subseteq\{\a_1,\cdots,\a_m\}$, then
\be \label{multiply-syzygy}
[\b_1\cdots\b_n][\c_1\cdots\c_n]=\sum_{(n-r,r)\vdash\C}[\b_1\cdots\b_r\C_{(1)}][\C_{(2)}\b_{r+1}\cdots\b_n].
\ee

Based on different choices of VW syzygies and multiple syzygies, White gave five different straightening algorithms. White also tested these five straightening algorithms and compared with the classical straightening algorithm. According to his results, all the five straightening algorithms have a higher efficiency than the classical straightening algorithm. And algorithm C and algorithm D are much better than the other three.

As we can see, the classical straightening algorithm and White's straightening algorithms can only deal with bracket polynomials of degree two, because all the syzygies have degree two in each step. Indeed, VW syzygy and multiple sygyzy are two different Gr\"{o}bner bases of the syzygy ideal of bracket algebra. So, the classical straightening algorithm and White's straightening algorithms are procedures of Gr\"{o}bner bases reduction.

\bn \cite{rota-straightening-formula}
The column order in bracket algebra is defined as: Let $[F], [G]$ be two straight bracket monomials, with degree $d$ and $l$ respectively, if
\begin{enumerate}
  \item $d<l$, or
  \item $d=l$, while there exist $i,j$, such that $F$ and $G$ contain the same first $(j-1)$ column. Denoting the $j$-th column of $F$ and $G$ are
  \[\begin{array}{c}
       \f_{1j} \\
       \f_{2j} \\
       \vdots \\
       \f_{dj} \\
     \end{array}~\textmd{and}~\begin{array}{c}
       \g_{1j} \\
       \g_{2j} \\
       \vdots \\
       \g_{lj} \\
     \end{array},\]
  then for any $k<i$, $\f_{kj}=\g_{kj}$, but $\f_{ij}\prec \g_{ij}$.
\end{enumerate}
Then we call $[F]$ has a low column order than $[G]$, denoted as $[F]\prec_c [G]$.
\en

\bn
The multiple order $\prec_m$ in bracket algebra is defined as: Let $[F], [G]$ be two bracket monomials, then we say $[F]\prec_m[G]$ if
\begin{enumerate}
  \item Under the negative column order, $[F^{\downarrow\downarrow}]\prec[G^{\downarrow\downarrow}]$, or
  \item $[F^{\downarrow\downarrow}]=[G^{\downarrow\downarrow}]$, and under the column order, $F\prec G$.
\end{enumerate}
\en

\bp
Multiple syzygies forms a Gr\"{o}bner basis of the syzygy ideal of bracket algebra under the multiple order.
\ep

\bo
Let $f$ be a bracket monomial of degree 2, and dimension $n$. For $1\leq j\leq n$, under the lexicographical orer, the smaller element in the $j$-th column of $f$, if appeared in the first row, then denoted as $\s_{1j}$; if appeared in the second row, then denoted as $\s_{2j}$. The larger element in the $j$-th column of $f$ will be denoted as $\m_j$. Assume that there are $r$ smaller elements in the first row. In (\ref{multiply-syzygy}), choosing
\begin{eqnarray*}
&& \c_1\c_2\cdots\c_n = \m_1\m_2\cdots\m_n; \\
&& \b_1\b_2\cdots\b_r = \s_{1j_1}\s_{1j_2}\cdots\s_{1j_r}; \\
&& \b_{r+1}\b_{r+2}\cdots\b_n = \s_{2j_{r+1}}\s_{2j_{r+2}}\cdots\s_{2j_n},
\end{eqnarray*}
we get a new expression of $f$. Under the multiple order, the leading term of this new expression is $f^{\downarrow\downarrow}$, which is straight. This means, multiple syzygies decreases the multiple order. Therefore, multiple syzygies forms a Gr\"{o}bner basis of the syzygy ideal of bracket algebra under the multiple order.
\hfill $\square$
\eo

\subsection{Rota's straightening algorithm}

Let $f$ be a monic bracket monomial, for any $1\leq i\leq m$, denote $\alpha_i$ as the number of $\a_i$ appeared in $f$, then
\be \label{content}
(\alpha)=(\alpha_1,\alpha_2,\cdots,\alpha_m)
\ee
is called the \emph{content} of $f$. Let $\{\r_1,\r_2,\cdots,$ $\r_n\}$ be $n$ vector variables with the property that $[\r_1\r_2\cdots\r_n]$ $=1$.

\bn
For any $l\geq 1$ and $\r_j,\a_i$, define the set polarization operators
\be \label{set polarization operator}
\mathcal{D}^l(\r_j,\a_i): \hspace{.5cm} \mathcal{B}[\{\a_u,\r_v\}]\longrightarrow\mathcal{B}[\{\a_u,\r_v\}]
\ee
as
\begin{enumerate}
  \item if $\alpha_i<l,$ then $\mathcal{D}^l(\r_j,\a_i)(f):=0;$
  \item if $\alpha_i\geq l,$ then $\mathcal{D}^l(\r_j,\a_i)(f):=\sum_{r=1}^{C_{\alpha_i}^l}f_r,$ where $f_1,f_2,\cdots,f_{C_{\alpha_i}^l}$ are all the $C_{\alpha_i}^l$ distinct bracket monomials obtained from $f$ by replacing each subset of $l$ vectors $\a_i$ by $l$ vectors $\r_j$.
\end{enumerate}
\en

\bx Let $f=\left[
                     \begin{array}{c}
                       \a_1\a_2\a_4 \\
                       \a_2\a_3\a_5 \\
                       \a_2\a_6\a_8 \\
                     \end{array}
                   \right]$, then
\[\mathcal{D}^2(\r_1,\a_2)(f)=\left[
                     \begin{array}{c}
                       \a_1\r_1\a_4 \\
                       \r_1\a_3\a_5 \\
                       \a_2\a_6\a_8 \\
                     \end{array}
                   \right]+\left[
                     \begin{array}{c}
                       \a_1\r_1\a_4 \\
                       \a_2\a_3\a_5 \\
                       \r_1\a_6\a_8 \\
                     \end{array}
                   \right]+\left[
                     \begin{array}{c}
                       \a_1\a_2\a_4 \\
                       \r_1\a_3\a_5 \\
                       \r_1\a_6\a_8 \\
                     \end{array}
                   \right].\]
\ex

The set polarization operators are commute with each other, which constitute the basic elements of Capelli operator.

\bn
Let $f$ be a monic bracket monomial, the Capelli operator generated by $f$ is
\be
\mathcal{C}(f):=\prod_{1\leq q\leq n}~\prod_{1\leq i\leq m}\mathcal{D}^{\eta_i(q)}(\r_q,\a_i),
\ee
where $\eta_i(q)$ equals to the the number of occurrences of $\a_i$ in the $q$-th column of $f$.
\en

\bx Let $f=\left[
           \begin{array}{c}
           \a_1\a_5\a_6 \\
           \a_2\a_3\a_7 \\
           \a_3\a_6\a_8 \\
           \end{array}
           \right]$, and $m=8,n=3$. Then
\begin{eqnarray*}
  && \eta_1(1)=\eta_2(1)=\eta_3(1)=1,~~~~~\eta_i(1)=0,~(i=4,5,6,7,8,9); \\
  && \eta_3(2)=\eta_5(2)=\eta_6(2)=1,~~~~~\eta_j(2)=0,~(j=1,2,4,7,8,9); \\
  && \eta_6(3)=\eta_7(3)=\eta_8(3)=1,~~~~~\eta_k(3)=0,~(k=1,2,3,4,5,6).
\end{eqnarray*}
So
\begin{eqnarray*}
  \mathcal{C}(f) &=& \mathcal{D}^1(\r_1,\a_1)\mathcal{D}^1(\r_1,\a_2)\mathcal{D}^1(\r_1,\a_3)\mathcal{D}^1(\r_2,\a_3)\mathcal{D}^1(\r_2,\a_5) \\
                 & & \mathcal{D}^1(\r_2,\a_6)\mathcal{D}^1(\r_3,\a_6)\mathcal{D}^1(\r_3,\a_8)\mathcal{D}^1(\r_3,\a_9).
\end{eqnarray*}
Operating this on $f$, we get $\mathcal{C}(f)(f)=\left[
                     \begin{array}{c}
                       \r_1\r_2\r_3 \\
                       \r_1\r_2\r_3 \\
                       \r_1\r_2\r_3 \\
                     \end{array}
                   \right]$. Let $g=\left[
                     \begin{array}{c}
                       \a_1\a_2\a_3 \\
                       \a_3\a_6\a_7 \\
                       \a_5\a_6\a_8 \\
                     \end{array}
                   \right]$, then $\mathcal{C}(f)(g)=0.$
\ex

\bl \cite{rota-straightening-formula}
Let $f$ and $g$ be two monic bracket monomials, with the same content, and $f\succ g$ under the negative column order, then
$\mathcal{C}(f)(f)\neq0$, and $\mathcal{C}(f)(g)=0$.
\el

More about Capelli operator can be found in \cite{clausen2}.

\bl \cite{rota-straightening-formula} \label{rota-straightening-formula}
Assume that $\{S_1,S_2,\cdots,S_N\}$ are the set of all straight bracket monomials with the same content as $f$. Under the negative column order $S_1\succ S_2\succ\cdots\succ S_N$, then there exist unique $\lambda_i~(i=1,2,\cdots,N)$, such that
\be \label{rota-straightening-algorithm}
f=\ds\sum \lambda_i S_i.
\ee
When operating $\mathcal{C}(S_i)$ on $f$, we will get a system of linear equations about $\lambda_1,\lambda_2,\cdots,\lambda_N$. The coefficient matrix is a lower triangular matrix and the diagonal elements are all nonzero.
\el

Therefore, we have the following Rota's straightening algorithm \cite{rota-straightening-formula}, \cite{desarmenien}:
\begin{algorithm}[H]
\caption{Straightening algorithm based on Capelli operator: \textsf{rota}}
\label{straightening-algorithm-based-on-Capelli-operator}
\begin{algorithmic}[1]
\REQUIRE
a homogeneous bracket polynomial $F$
\ENSURE
the straight expression of $F$
\STATE Find the set $\mathcal{S}=\{S_1,\cdots,S_N\}$ of all straight bracket monomials, where each $S_i$ has the same content as $F$, under the negative column order $S_i\succ S_{i+1}$. \\
\STATE Set $f:=F,~~~g:=0$.
\STATE For $i$ from 1 to $N$ do\\
       If $\mathcal{C}(S_i)(f)\neq0$, then  \\
       $g:=g+\mathcal{C}(S_i)(f)S_i,~~~f:=f-\mathcal{C}(S_i)(f)S_i.$ \\
       End if; End do.
\RETURN $g$.
\end{algorithmic}
\end{algorithm}

\bx  Let $f=\left[
   \begin{array}{c}
     \a_1\a_4\a_6 \\
     \a_2\a_3\a_5 \\
   \end{array}
 \right]$, then all straight bracket monomials with the same content as $f$ are:
\[f_1=\left[
   \begin{array}{c}
     \a_1\a_3\a_5 \\
     \a_2\a_4\a_6 \\
   \end{array}
 \right],~~~~f_2=\left[
   \begin{array}{c}
     \a_1\a_3\a_4 \\
     \a_2\a_5\a_6 \\
   \end{array}
 \right],~~~~f_3=\left[
   \begin{array}{c}
     \a_1\a_2\a_5 \\
     \a_3\a_4\a_6 \\
   \end{array}
 \right],~~~~f_4=\left[
   \begin{array}{c}
     \a_1\a_2\a_4 \\
     \a_3\a_5\a_6 \\
   \end{array}
 \right],~~~~f_5=\left[
   \begin{array}{c}
     \a_1\a_2\a_3 \\
     \a_4\a_5\a_6 \\
   \end{array}
 \right].\]
We assume that $f=\lambda_1f_1+\lambda_2f_2+\lambda_3f_3+\lambda_4f_4+\lambda_5f_5$.
By lemma \ref{rota-straightening-formula}, we have the following system of linear equations:
\[\left(
    \begin{array}{ccccc}
      1 & 0 & 0 & 0 & 0 \\
      0 & 1 & 0 & 0 & 0 \\
      0 & 0 & 1 & 0 & 0 \\
      0 & 0 & 0 & 1 & 0 \\
      1 & 0 & 0 & 0 & 1 \\
    \end{array}
  \right)\left(
           \begin{array}{c}
             \lambda_1 \\
             \lambda_2 \\
             \lambda_3 \\
             \lambda_4 \\
             \lambda_5 \\
           \end{array}
         \right)=\left(
           \begin{array}{c}
             \hfill 1 \\
             \hfill 0 \\
             -1 \\
             \hfill 0 \\
             \hfill 0 \\
           \end{array}
         \right),\]
with solution: $\lambda_1=1,~~\lambda_2=0,~~\lambda_3=-1,~~\lambda_4=0,~~\lambda_5=-1.$ Hence
\[f=\left[
   \begin{array}{c}
     \a_1\a_3\a_5 \\
     \a_2\a_4\a_6 \\
   \end{array}
 \right]-\left[
   \begin{array}{c}
     \a_1\a_2\a_3 \\
     \a_4\a_5\a_6 \\
   \end{array}
 \right]-\left[
   \begin{array}{c}
     \a_1\a_2\a_5 \\
     \a_3\a_4\a_6 \\
   \end{array}
 \right].\]
\ex

\section{Dual bracket and straightening algorithm based on it}
\setcounter{equation}{0}

In this section, we will present the concept of dual bracket, which is a new concept about Young tableau. Then some propositions about dual bracket will be given. Finally, we will show the new straightening algorithm in bracket algebra based on dual bracket, and describe its relationship with Rota's straightening algorithm.

\bn \cite{minc} (1). Let $A=(u_{ij})_{n\times n}$ be a matrix, the permanent of $A$ is defined as:
\[\emph{perm}(A)=\sum_{\sigma\in S_n}u_{1\sigma(1)}u_{2\sigma(2)}\cdots u_{n\sigma(n)}.\]
(2). Let $B=(v_{ij})_{r\times s}$ and $C=(w_{ij})_{r\times s}$ be two matrices with the same shape, the Hadamard product of these two matrices is defined as:
\[B\odot C=(v_{ij}w_{ij})_{r\times s}.\]
\en

\bn \label{dual bracket:def}
Suppose that $Y=\begin{array}{cccc}
                    \y_{11} & \y_{12} & \cdots & \y_{1n} \\
                    \y_{21} & \y_{22} & \cdots & \y_{2n} \\
                    \vdots  & \vdots  & \ddots & \vdots \\
                    \y_{d1} & \y_{d2} & \cdots & \y_{dn} \\
                  \end{array}$ is a Young tableau with degree $d$ and dimension $n$, where $\y_{ij}\in\{\a_1,\cdots,\a_m\}$. Let $H$ be the matrix
\be\label{matrix} \vspace{.1cm}
H=(\y_{11}\odot\y_{21}\cdots\odot\y_{d1},~~~~\y_{12}\odot\y_{22}\cdots\odot\y_{d1},~~~~\cdots,~~~~\y_{1n}\odot\y_{2n}\cdots\odot\y_{dn}).
\ee
Then define
\be\label{dual bracket:def formula}
\langle Y\rangle:=\emph{perm}(H), \hspace{1cm} \{Y\}:=\det(H),
\ee
where $\langle Y\rangle$ is called an even dual bracket of $Y$, $\{Y\}$ is called an odd dual bracket of $Y$, $d$ and $n$ are called the degree and dimension of dual bracket respectively.
\en

\bx (1). Assume that~$A=\left(
                         \begin{array}{cc}
                           x_1 & x_2 \\
                           x_3 & x_4 \\
                         \end{array}
                       \right)$, then~$\textmd{perm}(A)=x_1x_4+x_2x_3$.

(2). If~$B=\left(
                \begin{array}{cc}
                  y_1 & y_2 \\
                  y_3 & y_4 \\
                \end{array}
             \right)$, then~$A\odot B=\left(
                \begin{array}{cc}
                  x_1y_1 & x_2y_2 \\
                  x_3y_3 & x_4y_4 \\
                \end{array}
             \right)$.

(3). Suppose that~$n=d=2$, then
\[\left\{
     \begin{array}{cc}
       \a_1 & \a_2 \\
       \a_3 & \a_4 \\
     \end{array}
   \right\}=x_1x_3y_2y_4+x_2x_4y_1y_3, \hspace{1cm}
  \left\langle
     \begin{array}{cc}
       \a_1 & \a_2 \\
       \a_3 & \a_4 \\
     \end{array}
   \right\rangle=x_1x_3y_2y_4-x_2x_4y_1y_3.\]
\ex

The following is an easy proved property about dual bracket.

\bp \label{dual bracket:symmetry prop}
(1). An even dual bracket contains $n!$ terms, which is invariant if we permute two elements in a column, and also invariant under the permutation of two columns entirely.

(2). An odd dual bracket contains $n!$ terms, which is invariant if we permute two elements in a column, while anti-invariant under the permutation of two columns entirely.
\ep

\bn
Assume that~$Y=\C_1~\C_2~\cdots~\C_n$ is a Young tableau of degree $d$ and dimension $n$, where $\C_j$ is the $j$-th column of $Y$. Let $\sigma\in S_n$, then define
\be \label{Young tableau: column permutaion}
\sigma(T):= \C_{\sigma(1)}~\C_{\sigma(2)}~\cdots~\C_{\sigma(n)}.
\ee
\en

\bt (Multiplication Rules among Dual Brackets) \label{multiplication rules}
Let $C,D$ be two Young tableaux with the same dimension $n$, then
\be\label{multiplication rule formula}
\begin{tabular}{ll} \vspace{.2cm}
  \emph{(i1):}~~~$\langle C\rangle\langle D\rangle=\ds\sum_{\sigma\in S_n}\left\langle
                                                                     \begin{array}{c}
                                                                        C \\
                                                                       \sigma(D) \\
                                                                     \end{array}
                                                                   \right\rangle$,
&\hspace{1cm} \emph{(i2):}~~~$\langle C\rangle\langle D\rangle=\ds\sum_{\sigma\in S_n}\left\langle
                                                                     \begin{array}{c}
                                                                        \sigma(C) \\
                                                                        D \\
                                                                     \end{array}
                                                          \right\rangle$; \\ \vspace{.2cm}
  \emph{(ii1):}~~~$\langle C\rangle\{D\}=\ds\sum_{\sigma\in S_n}\emph{sign}(\sigma)\left\{
                                                               \begin{array}{c}
                                                                 C \\
                                                                 \sigma(D) \\
                                                               \end{array}
                                                           \right\}$,
&\hspace{1cm} \emph{(ii2):}~~~$\langle C\rangle\{D\}=\ds\sum_{\sigma\in S_n}\left\{
                                                  \begin{array}{c}
                                                     \sigma(C) \\
                                                     D \\
                                                  \end{array}
                                               \right\}$; \\
  \emph{(iii1):}~~~$\{C\}\{D\}=\ds\sum_{\sigma\in S_n}\emph{sign}(\sigma)\left\langle
                                                        \begin{array}{c}
                                                            C \\
                                                            \sigma(D) \\
                                                        \end{array}
                                                      \right\rangle$,
&\hspace{1cm} \emph{(iii2):}~~~$\{C\}\{D\}=\ds\sum_{\sigma\in S_n}\emph{sign}(\sigma)\left\langle
                                                         \begin{array}{c}
                                                             \sigma(C) \\
                                                             D \\
                                                         \end{array}
                                                      \right\rangle.$  \\
\end{tabular}
\ee
\et

\bo
We only give a proof of (i1), other proofs are similar.
Assume that the degree of $C$ and $D$ are $k$ and $l$ respectively,
\[C=\begin{array}{cccc}
        \c_{11} & \c_{12} & \cdots & \c_{1n} \\
        \c_{21} & \c_{22} & \cdots & \c_{2n} \\
        \vdots  & \vdots  & \ddots & \vdots  \\
        \c_{k1} & \c_{k2} & \cdots & \c_{kn} \\
      \end{array}, \hspace{2cm} D=\begin{array}{cccc}
        \d_{11} & \d_{12} & \cdots & \d_{1n} \\
        \d_{21} & \d_{22} & \cdots & \d_{2n} \\
        \vdots  & \vdots  & \ddots & \vdots  \\
        \d_{l1} & \d_{l2} & \cdots & \d_{ln} \\
      \end{array}.\]
Setting~$\c_{ij}=(x_{ij,1},x_{ij,2},\cdots,x_{ij,n})^T, \d_{uv}=(y_{uv,1},y_{uv,2},\cdots,y_{uv,n})^T$ be the coordinates, then
\begin{eqnarray*}
\ds\sum_{\sigma\in S_n}\left\langle
                            \begin{array}{c}
                                \sigma(C) \\
                                D \\
                            \end{array}
                         \right\rangle
&=& \sum_{\sigma\in S_n}\Bigg\{\sum_{\tau\in S_n}\Big( x_{1\sigma(1),\tau(1)} x_{2\sigma(1),\tau(1)} \cdots x_{k\sigma(1),\tau(1)}\cdots
                           x_{1\sigma(n),\tau(n)} x_{2\sigma(n),\tau(n)}~~~~~~~~~~\\
& & \cdots x_{k\sigma(n),\tau(n)}\Big)\times\Big(   y_{11,\tau(1)} y_{21,\tau(1)} \cdots y_{k1,\tau(1)}\cdots
                           y_{1n,\tau(n)} y_{2n,\tau(n)} \cdots y_{kn,\tau(n)}\Big)\Bigg\} \\
&=& \sum_{\tau\in S_n}\Bigg(\sum_{\sigma\in S_n}x_{1\sigma(1),\tau(1)} x_{2\sigma(1),\tau(1)} \cdots x_{k\sigma(1),\tau(1)}\cdots
                           x_{1\sigma(n),\tau(n)} x_{2\sigma(n),\tau(n)} \\
& & \cdots x_{k\sigma(n),\tau(n)}\Bigg)\times\Big(   y_{11,\tau(1)} y_{21,\tau(1)} \cdots y_{k1,\tau(1)}\cdots
                           y_{1n,\tau(n)} y_{2n,\tau(n)} \cdots y_{kn,\tau(n)}\Big) \\
&=& \sum_{\tau\in S_n}\Bigg(\sum_{\rho\in S_n}\Big( x_{11,\rho(1)} x_{21,\rho(1)} \cdots x_{k1,\rho(1)}\cdots
                           x_{1n,\rho(n)} x_{2n,\rho(n)} \cdots x_{kn,\tau(n)}\Big)\Bigg) \\
& & \hspace{1.5cm} \times~~\Big(   y_{11,\tau(1)} y_{21,\tau(1)} \cdots y_{k1,\tau(1)}\cdots
                           y_{1n,\tau(n)} y_{2n,\tau(n)} \cdots y_{kn,\tau(n)}\Big) \\
&=& ~~~~~\Big(\sum_{\rho\in S_n}x_{11,\rho(1)} x_{21,\rho(1)} \cdots x_{k1,\rho(1)}\cdots
                           x_{1n,\rho(n)} x_{2n,\rho(n)} \cdots x_{kn,\tau(n)}\Big) \\
& & \times~~\Big(\sum_{\tau\in S_n}y_{11,\tau(1)} y_{21,\tau(1)} \cdots y_{k1,\tau(1)}\cdots
                           y_{1n,\tau(n)} y_{2n,\tau(n)} \cdots y_{kn,\tau(n)}\Big) \\
&=& \langle C\rangle\langle D\rangle.
\end{eqnarray*}
The theorem has been proved.
\hfill $\square$
\eo

This theorem depicts a reason why we call $\langle T\rangle$ even dual bracket, $\{T\}$ odd dual bracket: the multiplication rules among dual brackets are similar to the addition rules among even and odd integral numbers.

\bn
Assume that~$Y=\begin{array}{c}
                   \R_1 \\
                   \R_2 \\
                   \vdots \\
                   \R_d \\
                 \end{array}$ is a Young tableau of degree $d$ and dimension $n$, where $\R_i$ is the $i$-th row of $Y$. Let $\sigma_1, \sigma_2,\cdots,\sigma_d\in S_n$. Define
\be \label{Young tableau: single column permutaion}
Y_{\sigma_1,\sigma_2,\ldots,\sigma_d}:=\begin{array}{c}
                                 \sigma_1(\R_1) \\
                                 \sigma_2(\R_2) \\
                                 \vdots \\
                                 \sigma_d(\R_d) \\
                               \end{array}.
\ee
\en

\bn \label{db:def}
Let~$Y$ be a Young tableau with degree $d$ and dimension $n$. Define
\be \label{dual tableau polynomial}
DT(Y):=\sum_{\sigma_2,\cdots,\sigma_d\in S_n}\emph{sign}(\sigma_2)\cdots \emph{sign}(\sigma_d)~Y_{\emph{id},\sigma_2,\ldots,\sigma_d},
\ee
which is called the dual tableau polynomial of $Y$. Moreover, define
\be
DB([Y])=[DT(Y)],
\ee
which is called the dual bracket polynomial of $[Y]$.
\en

\bn \label{sb:def}
Let~$Y$ be a Young tableau with degree $d$ and dimension $n$. Setting
\be
\mathcal{S}=\Big\{(\tau_2,\cdots,\tau_d)\in S_n^{d-1}~\Big|~\mathfrak{N}_C(Y_{\emph{id},\tau_2,\ldots,\tau_d})~\textmd{is~straight}\Big\},
\ee
where~$\mathfrak{N}_C(Y_{\emph{id},\tau_2,\ldots,\tau_d})$ is the column normalized tableau of $Y_{\emph{id},\tau_2,\ldots,\tau_d}$. Define
\be \label{sb}
ST(Y):=\sum_{(\tau_2,\cdots,\tau_d)\in \mathcal{S}} \emph{sign}(\tau_2)\cdots \emph{sign}(\tau_d)~\mathfrak{N}_C(Y_{\emph{id},\tau_2,\ldots,\tau_d}),
\ee
which is called the special tableau polynomial of $Y$. Also define
\be \label{special bracket polynomial}
SB([Y]):=[ST(Y)],
\ee
which is called the special bracket polynomial of the bracket monomial $[Y]$.
\en

\bx Let~$Y=\begin{array}{ccc}
     \a_1 & \a_3 & \a_5 \\
     \a_2 & \a_4 & \a_6 \\
   \end{array}, d=2,n=3$, then
\begin{eqnarray*}
  \mathfrak{N}_C(DT(Y)) &=&
   \begin{array}{ccc}
     \a_1 & \a_3 & \a_5 \\
     \a_2 & \a_4 & \a_6 \\
   \end{array}-
   \begin{array}{ccc}
     \a_1 & \a_3 & \a_4 \\
     \a_2 & \a_6 & \a_5 \\
   \end{array}-
   \begin{array}{ccc}
     \a_1 & \a_2 & \a_5 \\
     \a_4 & \a_3 & \a_6 \\
   \end{array}+
   \begin{array}{ccc}
     \a_1 & \a_2 & \a_3 \\
     \a_4 & \a_5 & \a_6 \\
   \end{array}+
   \begin{array}{ccc}
     \a_1 & \a_2 & \a_4 \\
     \a_6 & \a_3 & \a_5 \\
   \end{array}-
   \begin{array}{ccc}
     \a_1 & \a_2 & \a_3 \\
     \a_6 & \a_5 & \a_4 \\
   \end{array}, \\
ST(Y) &=& \begin{array}{ccc}
     \a_1 & \a_3 & \a_5 \\
     \a_2 & \a_4 & \a_6 \\
   \end{array}+\begin{array}{ccc}
     \a_1 & \a_2 & \a_3 \\
     \a_4 & \a_5 & \a_6 \\
   \end{array}.
\end{eqnarray*}
\ex
This example shows that $ST(Y)$ has fewer terms than $\mathfrak{N}_C(DT(Y))$.

\bt \label{expanding-formulas}
Let~$Y$ be a Young tableau with degree $d$ and dimension $n$. Then 
\be \label{dual bracket:expanding formula}
[Y]^c=\left\{
  \begin{array}{ll} \vspace{.2cm}
    \langle DT(Y) \rangle, & \hbox{~~~if~$d$ is even;} \\
    \{DT(Y)\},             & \hbox{~~~if~$d$ is odd.}
  \end{array}
\right.
\ee
\et

\bo By induction on $d$. If $d=1$, then $[Y]^c=\{Y\}$, and thus~(\ref{dual bracket:expanding formula}) holds apparently. Assume that (\ref{dual bracket:expanding formula}) holds for $d=r$. We next prove (\ref{dual bracket:expanding formula}) holds for $d=r+1$. Denote $Y=S \circ_C R$, where $S$ is a Young tableau constituted by the first $r$ rows of $Y$, and $R$ is the remaining Young tableau in $Y$ with degree 1.

If~$r$ is even, then
\begin{eqnarray*}
[Y]^c &=& \sum_{\sigma_2,\cdots,\sigma_r\in S_n}\textmd{sign}(\sigma_2)\cdots \textmd{sign}(\sigma_r)\langle S_{\textmd{id},\sigma_2,\ldots,\sigma_r}\rangle \{R\} \hspace{4.9cm} (\textmd{By induction}) \\
      &=& \sum_{\sigma_2,\cdots,\sigma_r,\sigma_{r+1}\in S_n}\textmd{sign}(\sigma_2)\cdots \textmd{sign}(\sigma_r)\textmd{sign}(\sigma_{r+1})\left\{
                                 \begin{array}{c}
                                   S_{\textmd{id},\sigma_2,\ldots,\sigma_r} \\
                                   \sigma_{r+1}(R) \\
                                 \end{array}
                               \right\} \hspace{2.3cm} (\textmd{By (ii1) of}~(\ref{multiplication rule formula})) \\
      &=& \ds\sum_{\sigma_2,\cdots,\sigma_{r+1}\in S_n}\textmd{sign}(\sigma_2)\cdots \textmd{sign}(\sigma_{r+1})\{Y_{\textmd{id},\sigma_2,\ldots,\sigma_{r+1}}\}.
\end{eqnarray*}

If $r$ is odd, then
\begin{eqnarray*}
[Y]^c &=& \sum_{\sigma_2,\cdots,\sigma_r\in S_n}\textmd{sign}(\sigma_2)\cdots \textmd{sign}(\sigma_r)\{S_{\textmd{id},\sigma_2,\ldots,\sigma_r}\}\{R\} \hspace{4.8cm} (\textmd{By induction}) \\
      &=& \sum_{\sigma_2,\cdots,\sigma_r,\sigma_{r+1}\in S_n}\textmd{sign}(\sigma_2)\cdots \textmd{sign}(\sigma_r)\textmd{sign}(\sigma_{r+1})\left\langle
                                 \begin{array}{c}
                                   S_{\textmd{id},\sigma_2,\ldots,\sigma_r} \\
                                   \sigma_{r+1}(R) \\
                                 \end{array}
                               \right\rangle \hspace{2.2cm} (\textmd{By (iii1) of}~(\ref{multiplication rule formula})) \\
      &=& \ds\sum_{\sigma_2,\cdots,\sigma_{r+1}\in S_n}\textmd{sign}(\sigma_2)\cdots \textmd{sign}(\sigma_{r+1})\langle Y_{\textmd{id},\sigma_2,\ldots,\sigma_{r+1}}\rangle.
\end{eqnarray*}
Therefore, (\ref{dual bracket:expanding formula}) is right.
\hfill $\square$
\eo

\bl \label{leadingterm of SB}
Let~$Y$ be a Young tableau, then under the negative column order, the leading term of $SB([Y])$ is~$[Y]^{\downarrow\downarrow}$. So for any two Young tableaux $X,Y$, $SB([X])\succ SB([Y])$ under the negative column order if and only if $lt([X]^c)\succ lt([Y]^c)$ under the negative basis order.
\el

\bo We only need to prove that under the negative coordinate order, the leading term $lt(SB([Y])^c)$ of the coordinate expanding expression $SB([Y])^c$ of $SB([Y])$ equals to the leading term $lt(([Y]^{\downarrow\downarrow})^c)$ of the coordinate expanding expression $([Y]^{\downarrow\downarrow })^c$ of $[Y]^{\downarrow\downarrow}$. This is verified by the following three facts:

(a). $[Y]^{\downarrow\downarrow} $ is a term of $SB([Y])$,

(b). $lt([Y]^c)=lt(([Y]^{\downarrow\downarrow})^c)$,

(c). By~(\ref{sb}) we know that $lt(SB([Y])^c)$ equals to $lt([\mathfrak{N}_C(Y_{\textmd{id},\tau_2,\ldots,\tau_d})]^c)$ for some $\tau_2,\cdots,\tau_d\in S_n$.
In addition, by~(\ref{dual tableau polynomial}) we know that $lt([\mathfrak{N}_C(Y_{\textmd{id},\tau_2,\ldots,\tau_d})]^c)$ equals to either a term of $\langle DT(Y\rangle$ or a term of
$\{ DT(Y)\}$, this depends on the degree. Therefore, from (\ref{dual bracket:expanding formula}), we have $lt([\mathfrak{N}_C(Y_{\textmd{id},\tau_2,\ldots,\tau_d})]^c)$ equals to a term of the coordinate expanding expression of $[Y]$.
\hfill $\square$
\eo

Based on the relationship between bracket algebra and its coordinate polynomial ring, i.e., lemma \ref{relationship between bracket algebra and its coordinate ring}, we can easily get the following straightening algorithm about bracket polynomials.

\begin{algorithm}[H]
\caption{Straightening algorithm based on coordinate ring}
\label{straightening algorithm based on coordinate ring}
\begin{algorithmic}[1]
\REQUIRE
A homogeneous bracket polynomial $F=\sum_{i=1}^M\alpha_i[X_i]$
\ENSURE
The straightening expression of $F$
\STATE For $i$ from 1 to $M$ do \\
       compute the coordinate expanding expression $[X_i]^c$ of $[X_i]$ \\
       end do
\STATE Set $p:=\sum \alpha_i[X_i]^c,~~~g:=0$.
\STATE While~$p\neq0$ do \\
       $M:=lt(p),~~~p:=p-[Tableau(M)]^c,~~~g:=g+[Tableau(M)].$ \\
       End do.
\RETURN $g$.
\end{algorithmic}
\end{algorithm}

Theorem \ref{expanding-formulas} together with lemma \ref{leadingterm of SB} imply that we can also study bracket polynomials from dual bracket, and this can get the same results as from coordinate rind. Based on algorithm \ref{straightening algorithm based on coordinate ring}, we can get the following straightening algorithm about bracket polynomials.

\begin{algorithm}[H]
\caption{Straightening algorithm based on coordinate ring: \textsf{db}}
\label{straightening algorithm based on dual bracket}
\begin{algorithmic}[1]
\REQUIRE
A homogeneous bracket polynomial $F=\sum_{i=1}^M\alpha_i[X_i]$
\ENSURE
The straightening expression of $F$
\STATE For $i$ from 1 to $M$ do \\
       compute $SB([X_i])$ \\
       end do
\STATE Set $q:=\sum \alpha_iSB([X_i]),~~~g:=0$.
\STATE While~$p\neq0$ do \\
       $L:=LC(q)LM(q),~~~q:=q-LC(q)SB(LM(q)),~~~g:=g+L$. \\
       End do.
\RETURN $g$.
\end{algorithmic}
\end{algorithm}

\bp \label{dual algorithm validity lemma}
Algorithm~\ref{straightening algorithm based on dual bracket} terminates in finite steps and return the straight expression of $F$.
\ep

\bo Due to the similar procedure of algorithm \ref{straightening algorithm based on coordinate ring} and algorithm \ref{straightening algorithm based on dual bracket}. By lemma \ref{leadingterm of SB}, for any Young tableau $R$, $lt(SB([R]]^c)=lt([R]^c)$. So in algorithm \ref{straightening algorithm based on dual bracket}, the leading term of the coordinate expanding expression of $LT(q)=LC(q)LM(q)$ is $lt(LT(q)^c)$, which equals to $lt(p)$ in algorithm \ref{straightening algorithm based on coordinate ring}. Since algorithm \ref{straightening algorithm based on coordinate ring} terminates in finite steps and return the straight expression of $F$, and so is algorithm \ref{straightening algorithm based on dual bracket}.
\hfill $\square$
\eo

\bx Let $f=\left[\begin{array}{c}
           \a_1  \a_4  \a_6 \\
           \a_2  \a_3  \a_5 \\
         \end{array}\right]$, then

\textbf{Step 1}, by lemma~\ref{leading term:bracket monomial}, the leading term of the straight expression of $f$ under the negative column order is
~$f_1=f^{\downarrow\downarrow}=\left[\begin{array}{c}
           \a_1  \a_3  \a_5 \\
           \a_2  \a_4  \a_6 \\
         \end{array}\right]$. Then we need to calculate $SB(f)-SB(f_1)$. By definition
\[SB(f)=\left[\begin{array}{c}
           \a_1  \a_3  \a_5 \\
           \a_2  \a_4  \a_6 \\
         \end{array}\right]-
         \left[\begin{array}{c}
           \a_1  \a_2  \a_5 \\
           \a_3  \a_4  \a_6 \\
         \end{array}\right];~~~SB(f_1)=
         \left[\begin{array}{c}
           \a_1  \a_3  \a_5 \\
           \a_2  \a_4  \a_6 \\
         \end{array}\right]+
         \left[\begin{array}{c}
           \a_1  \a_2  \a_3 \\
           \a_4  \a_5  \a_6 \\
         \end{array}\right].\]
So
\[SB(f)-SB(f_1)=-\left[\begin{array}{c}
           \a_1  \a_2  \a_5 \\
           \a_3  \a_4  \a_6 \\
         \end{array}\right]-
         \left[\begin{array}{c}
           \a_1  \a_2  \a_3 \\
           \a_4  \a_5  \a_6 \\
         \end{array}\right].\]
Hence the leading term of $SB(f)-SB(f_1)$ under  the negative column order is $f_2=-
         \left[\begin{array}{c}
           \a_1  \a_2  \a_5 \\
           \a_3  \a_4  \a_6 \\
         \end{array}\right]$.

\textbf{Step 2}, compute $SB(f_2)=f_2$, so $SB(f)-SB(f_1)-SB(f_2)=-
         \left[\begin{array}{c}
           \a_1  \a_2  \a_3 \\
           \a_4  \a_5  \a_6 \\
         \end{array}\right]$, whose leading term under  the negative column order is $f_3=-
         \left[\begin{array}{c}
           \a_1  \a_2  \a_3 \\
           \a_4  \a_5  \a_6 \\
         \end{array}\right]$.

\textbf{Step 3}, compute $SB(f_3)=f_3$, so $SB(f)-SB(f_1)-SB(f_2)-SB(f_3)=0$, algorithm terminates. Finally we get
\[f=f_1+f_2+f_3=
   \left[\begin{array}{c}
     \a_1  \a_3  \a_5 \\
     \a_2  \a_4  \a_6 \\
   \end{array}\right]-
   \left[\begin{array}{c}
     \a_1  \a_2  \a_3 \\
     \a_4  \a_5  \a_6 \\
   \end{array}\right]-
   \left[\begin{array}{c}
     \a_1  \a_2  \a_5 \\
     \a_3  \a_4  \a_6 \\
   \end{array}\right].\]
\ex

The following proposition describes the relationship between dual bracket and Capelli operator, and then describes the relationship between algorithm \ref{straightening algorithm based on dual bracket} and Rota's straightening algorithm.

\bp Let $f,g$ be two bracket monomials with the same content. Assume that $\mathcal{C}(f)(g)\neq0$, and $DB(g)=\sum\lambda_i g_i$. Then there exists an unique $i$, such that $f=g_i$ and $\mathcal{C}(f)(g)=\lambda_i$.
\ep

\bo Let $f=[F]=[\f_{ij}]$ and $g=[G]=[\g_{ij}]$, where $\f_{ij},\g_{ij}\in\{\a_1,\a_2,\cdots,\a_m\}$. Assume that they have degree $d$ and dimension $n$.

\textbf{Case 1, multilinear case.} Since $f,g$ have the same content, so for any $1\leq i\leq d$ and $1\leq j\leq n$, there exist unique $1\leq u_{ij}\leq d$ and $1\leq v_{ij}\leq n$, such that $\g_{ij}=\f_{u_{ij}v_{ij}}$. So
\[g=\left[
       \begin{array}{cccc}
         \f_{u_{11}v_{11}} & \f_{u_{12}v_{12}} & \cdots & \f_{u_{1n}v_{1n}} \\
         \f_{u_{21}v_{21}} & \f_{u_{22}v_{22}} & \cdots & \f_{u_{2n}v_{2n}} \\
         \vdots            & \vdots            & \ddots & \vdots            \\
         \f_{u_{d1}v_{d1}} & \f_{u_{d2}v_{d2}} & \cdots & \f_{u_{dn}v_{dn}} \\
       \end{array}
     \right]~\textmd{and}~\mathcal{C}(f)(g)=\left[
                      \begin{array}{cccc}
                        \r_{v_{11}} & \r_{v_{12}} & \cdots & \r_{v_{1n}} \\
                        \r_{v_{21}} & \r_{v_{22}} & \cdots & \r_{v_{2n}} \\
                        \vdots      & \vdots      & \ddots & \vdots      \\
                        \r_{v_{d1}} & \r_{v_{d2}} & \cdots & \r_{v_{dn}} \\
                      \end{array}
                    \right].\]

Since $\mathcal{C}(f)(g)\neq0$, so for any $1\leq i \leq d$, there is an unique $\sigma_i\in S_n$, such that
\[(1,2,\cdots,n)=(v_{i\sigma_i(1)},v_{i\sigma_i(2)},\cdots,v_{i\sigma_i(n)}).\]
Then we have $\mathcal{C}(f)(g)=\textmd{sign}(\sigma_1)\textmd{sign}(\sigma_2)\cdots\textmd{sign}(\sigma_d)$. Therefore,
$F^{\downarrow\downarrow}=G_{\sigma_1,\sigma_2,\cdots,\sigma_d}^{\downarrow\downarrow}.$ Note that
\[G_{\sigma_1,\sigma_2,\cdots,\sigma_d}=\sigma_1(G_{\textmd{id},\sigma_1^{-1}\sigma_2,\cdots,\sigma_1^{-1}\sigma_d}),\]
and $G_{\textmd{id},\sigma_1^{-1}\sigma_2,\cdots,\sigma_1^{-1}\sigma_d}$ is a term of $DB(g)$.
Because of proposition \ref{dual bracket:symmetry prop}, if $d$ is even, then the coefficient is
\[\textmd{sign}(\sigma_1)^{d-1}\textmd{sign}(\sigma_2)\cdots\textmd{sign}(\sigma_d)=\textmd{sign}(\sigma_1)\textmd{sign}(\sigma_2)\cdots\textmd{sign}(\sigma_d)
=\mathcal{C}(f)(g);\]
and if $d$ is odd, then the coefficient is
\[\textmd{sign}(\sigma_1)^d\textmd{sign}(\sigma_2)\cdots\textmd{sign}(\sigma_d)=\textmd{sign}(\sigma_1)\textmd{sign}(\sigma_2)\cdots\textmd{sign}(\sigma_d)
=\mathcal{C}(f)(g).\]
Therefore, there is a $i$, such that $f=g_i$ and $\lambda_i=\mathcal{C}(f)(g)$.

\textbf{Case 2, general case.} Since $f,g$ have the same content, so for any $1\leq i\leq d$ and $1\leq j\leq n$, there exists an unique $k_{ij}$, such that for any $1\leq x_{ij}\leq k_{ij}$, we have $g_{ij}=f_{u_{ijx_{ij}}v_{ijx_{ij}}}$. So
\begin{eqnarray*}
g &=& \left[
       \begin{array}{cccc}
         \f_{u_{11x_{11}}v_{11x_{11}}}  & \f_{u_{12x_{12}}v_{12x_{12}}} & \cdots & \f_{u_{1nx_{1n}}v_{1nx_{1n}}}  \\
         \f_{u_{21x_{21}}v_{21x_{21}}} & \f_{u_{22x_{22}}v_{22x_{22}}} & \cdots & \f_{u_{2nx_{2n}}v_{2nx_{2n}}}  \\
         \vdots                                       & \vdots                                       & \ddots & \vdots            \\
         \f_{u_{d1x_{d1}}v_{d1x_{d1}}} & \f_{u_{d2x_{d2}}v_{d2x_{d2}}} & \cdots & \f_{u_{dnx_{dn}}v_{dnx_{dn}}}  \\
       \end{array}
     \right], \\
 \mathcal{C}(f)(g) &=& \sum_{1\leq i\leq d}~\sum_{1\leq j\leq n}~\sum_{1\leq x_{ij}\leq k_{ij}}\left[
                      \begin{array}{cccc}
                        \r_{v_{11x_{11}}} & \r_{v_{12x_{12}}} & \cdots & \r_{v_{1nx_{1n}}} \\
                        \r_{v_{21x_{21}}} & \r_{v_{22x_{22}}} & \cdots & \r_{v_{2nx_{2n}}} \\
                        \vdots            & \vdots            & \ddots & \vdots      \\
                        \r_{v_{d1x_{d1}}} & \r_{v_{d2x_{d2}}} & \cdots & \r_{v_{dnx_{dn}}} \\
                      \end{array}
                    \right].
\end{eqnarray*}
Since  $\mathcal{C}(f)(g)\neq 0$, so there exist $x_{ij}$, such that
\[s:=\left[
                      \begin{array}{cccc}
                        \r_{v_{11x_{11}}} & \r_{v_{12x_{12}}} & \cdots & \r_{v_{1nx_{1n}}} \\
                        \r_{v_{21x_{21}}} & \r_{v_{22x_{22}}} & \cdots & \r_{v_{2nx_{2n}}} \\
                        \vdots            & \vdots            & \ddots & \vdots      \\
                        \r_{v_{d1x_{d1}}} & \r_{v_{d2x_{d2}}} & \cdots & \r_{v_{dnx_{dn}}} \\
                      \end{array}
                    \right]\neq 0.\]
So for any $1\leq i \leq d$, there is an unique $\sigma_i\in S_n$, such that
\[(1,2,\cdots,n)=(v_{i\sigma_i(1)x_{i\sigma_i(1)}},v_{i\sigma_i(2)x_{i\sigma_i(2)}},\cdots,v_{i\sigma_i(n)x_{i\sigma_i(n)}}).\]
Then we have $s=\textmd{sign}(\sigma_1)\textmd{sign}(\sigma_2)\cdots\textmd{sign}(\sigma_d)$. Therefore,
$F^{\downarrow\downarrow}=G_{\sigma_1,\sigma_2,\cdots,\sigma_d}^{\downarrow\downarrow}.$ Note that
\[G_{\sigma_1,\sigma_2,\cdots,\sigma_d}=\sigma_1(G_{\textmd{id},\sigma_1^{-1}\sigma_2,\cdots,\sigma_1^{-1}\sigma_d}),\]
and $G_{\textmd{id},\sigma_1^{-1}\sigma_2,\cdots,\sigma_1^{-1}\sigma_d}$ is a term of $DB(g)$.
Similarly, because of proposition \ref{dual bracket:symmetry prop}, if $d$ is even, then the coefficient is
\[\textmd{sign}(\sigma_1)^{d-1}\textmd{sign}(\sigma_2)\cdots\textmd{sign}(\sigma_d)=\textmd{sign}(\sigma_1)\textmd{sign}(\sigma_2)\cdots\textmd{sign}(\sigma_d)=s;\]
and if $d$ is odd, then the coefficient is
\[\textmd{sign}(\sigma_1)^d\textmd{sign}(\sigma_2)\cdots\textmd{sign}(\sigma_d)=\textmd{sign}(\sigma_1)\textmd{sign}(\sigma_2)\cdots\textmd{sign}(\sigma_d)=s.\]
Therefore, adding them all we also get that there exists a $i$, such that $f=g_i$ and $\lambda_i=\mathcal{C}(f)(g)$. Thus, this proposition is right.
\hfill $\square$
\eo

As a corollary, we have
\bc
Let $f$ be a bracket monomial of degree $d$. Let $\{S_1,S_2,\cdots,S_N\}$ be the set of all straight bracket monomial of degree $d$, with the property that $S_i\succ S_{i+1}$ under the negative column order. Then
\be
SB(f)=\sum_{i=1}^N \mathcal{C}(S_i)(f)S_i.
\ee
\ec
%

The essential step in algorithm \ref{straightening algorithm based on dual bracket} relies on the computation of $SB([X])=[ST(X)]$, where $X$ is a Young tableau. In the following, we first present an algorithm to compute $ST(X)$ in the multi-linear case, i.e., all the vector variables in $X$ are different.

Let $X$ be a Young tableau of degree $d$ and dimension $n$. If $n=1$, then it's obvious $ST(X)=X^{\downarrow\downarrow}$. We now assume that $n>1$. Denote all the vector variables appeared in $X$ are $\{\d_1,\d_2,\cdots,\d_{dn}\}$, and $\d_i\prec \d_{i+1},~(1\leq i\leq dn-1)$ lexicographically. For any $1\leq s\leq d$, set
\be
\mathcal{U}_s=\{\d_1,\d_2,\cdots,\d_{(s-1)n+1}\}.
\ee
Define the \emph{first column set} of $ST(X)$ as
\be\ba{lll} \label{first-column}
FC(X):=\Big\{\begin{array}{c}
                  \u_1 \\
                  \u_2 \\
                  \vdots \\
                  \u_d \\
                \end{array} &\Big|& \textmd{for~all}~1\leq s\leq d,~\u_s\in\mathcal{U}_s,~\u_s\succ \u_{s-1}, ~\textmd{in}~T, \u_1,\u_2,\cdots,\u_d~\textmd{do~not~lie}~~~\\
                && \hfill \textmd{in~the~same~row}\Big\}.
\ea\ee

\bn
Let $V$ be a Young tableau of degree $d$ dimension 1, $W$ be a Young tableau of degree $d$ dimension $k$. For any $1\leq i\leq d$, assume that the $i$-th element of $V$ lies in the $p_i$-th position in the $i$-th row of $W$.

(1). In the $i$-th row of $W$, delete the $p_i$-th element, the result is denoted as $W\backslash V$, called delete $V$ in $W$ by rows.

(2). For any $1\leq i\leq d$, move all the elements after the $p_i$-th position forward in $W\backslash V$, the result is denoted as $DCR(W,V)$, called delete and compress $V$ in $W$ by rows.
\en

Denote all elements in $FC(X)$ as $C_1,C_2,\cdots,C_\alpha$. For any $1\leq i\leq \alpha$, denote $D_i=\mathfrak{N}_R(DCR(X,C_i))$, which is a Young tableau of dimension $(n-1)$. By induction, for any $1\leq i\leq \alpha$, we assume that $ST(D_i)$ has been computed, denoted as $ST(D_i)=\sum_{j=1}^{\xi_i}\lambda_{ij}E_{ij}$. For $1\leq i\leq \alpha$, define
\be \label{sum-choice}
\Lambda_i=\Big\{1\leq j\leq \xi_i~~\Big|~~C_i\circ_R E_{ij}~\textmd{is~straight}\Big\}.
\ee
For any $1\leq i\leq \alpha, 1\leq j\leq d$, assume the $j$-th element of $C_i$ lies in the $q_{ij}$-th column of $X$, define
\be \label{column-sign}
\epsilon_i=\sum_{j=1}^d (q_{ij}-1).
\ee
Denoting the first row of $X$ as $\t_1~\t_2~\cdots~\t_n$. For any $1\leq k\leq n$, let $\t_k$ appears in the $\omega_{ijk}$-th column of $C_i\circ_R E_{ij}$,
define the permutation $\pi_{ij}\in S_n$ as
\be \label{permutation}
\pi_{ij}(k)=\omega_{ijk},~(\forall~k=1,2,\cdots,n).
\ee

\bl \label{perm-lemma}
Let $\Omega=\{\omega_1,\omega_2,\cdots,\omega_k\}\subseteq \{1,2,\cdots,n\}$ be a subset, define the shift permutation of $\Omega$ as
\be
s_\Omega:~~\big(\omega_1,\omega_2,\cdots,\omega_{k-1},\omega_k\big)~~\mapsto~~\big(\omega_k,\omega_1,\omega_2,\cdots,\omega_{k-1}\big),
\ee
By keeping the elements outside of $\Omega$ invariant, $s_\Omega$ can be extended as a permutation of $\{1,2,\cdots,n\}$, which is also denoted as $s_\Omega$. Then any permutation $\rho\in S_n$ can be factored as the composition of $n$ shift permutation.
\el

\bo
For any $1\leq k\leq n$, set $\Phi_k=\{1,2,\cdots,\rho(k)\}$, $\Omega_k=\Phi_k-\Phi_k\cap\{\rho(1),\cdots,\rho(k-1)\}$. Then it's easy to verify that
\begin{eqnarray*}
 (1,2,\cdots,n)  & \xrightarrow[]{s_{\Omega_1}} & (\rho(1),1,2,\cdots,\check{\rho(1)},\cdots,n) \\
                 & \xrightarrow[]{s_{\Omega_2}} & (\rho(1),\rho(2),1,2,\cdots,\check{\rho(1)},\cdots,\check{\rho(2)},\cdots,n) \\
                 &                              & \cdots\cdots\cdots \\
                 & \xrightarrow[]{s_{\Omega_n}} & (\rho(1),\rho(2),\cdots,\rho(n))
\end{eqnarray*}
So $\rho=s_{\Omega_n}\circ s_{\Omega_{n-1}}\circ\cdots\circ s_{\Omega_1}$, and $\textmd{sign}(\rho)=(-1)^{\#(\Omega_1)-1}\cdots(-1)^{\#(\Omega_n)-1}$.
\hfill $\square$
\eo

\bp \label{sb-prop}
Notation as above, then
\be \label{sb-alg}
ST(X)=\sum_{i=1}^\alpha~\sum_{j\in\Lambda_i}\emph{sign}(\pi_{ij})^d(-1)^{\epsilon_i}\lambda_{ij}~~C_i\circ_R E_{ij}.
\ee
\ep

\bo Setting $ST(X)=\sum\zeta_k~P_k\circ_R Q_k$, where the dimension of $P_k$ and $Q_k$ are 1 and $(n-q)$ respectively.
By definition \ref{sb:def}, there exist $\rho_2,\cdots,\rho_d\in S_n$, such that
\[P_k\circ_R Q_k=\mathfrak{N}_C(X_{\rho_2,\ldots,\rho_d})=\begin{array}{cccc}
                                                          \r_{11} & \r_{12} & \cdots & \r_{1n} \\
                                                          \r_{21} & \r_{22} & \cdots & \r_{2n} \\
                                                          \vdots  & \vdots  & \ddots & \vdots  \\
                                                          \r_{d1} & \r_{d2} & \cdots & \r_{dn} \\
                                                        \end{array}.\]
For any $1\leq s\leq d$, in $P_k\circ_R Q_k$ the elements has a low order than $\r_{s1}$ must appears in the first $(s-1)$ rows, which contains $(s-1)n$ elements totally, so~$\r_{s1}\in\mathcal{U}_s$. By the definition of $ST(X)$, under the lexicographical order $\r_{s1}\succ \r_{(s-1)1}$ and $\r_{11},\r_{21},\cdots,\r_{d1}$ can not appear in the same row of $X$, so $P_k\in FC(T)$, which means the first column of each elements of $ST(X)$ comes from $FC(X)$. Similarly, the first column of each $Q_k$ also comes from some $FC(D_{l_k})$. Therefore, there exist $i,j$, such that $P_k\circ_R Q_k=C_i\circ_R E_{ij}$.

By lemma \ref{perm-lemma}, if set $\tau_2,\cdots,\tau_d\in S_n$, then $\epsilon_i\lambda_{ij}=\textmd{sign}(\tau_2)\cdots \textmd{sign}(\tau_d)$.

Finally, by lemma \ref{dual bracket:symmetry prop}, in (\ref{sb-alg}), $\textmd{sign}(\pi_{ij})^d$ is necessary.
\hfill $\square$
\eo

By the above analysis and proposition \ref{sb-prop}, we have the following algorithm to compute $ST(X)$ in the multi-linear case.
\begin{algorithm}[H]
\caption{Compute $ST(X)$ in the multi-linear case}
\label{a algorithm of compututation SB in the multilinear case}
\begin{algorithmic}[1]
\REQUIRE
A Young tableau $X$ of degree $d$ and dimension $n$
\ENSURE
$ST(X)$
\STATE Enumerate all the vector variables $\{\d_1,\cdots,\d_{dn}\}$ in $X$, where $\d_i\succ \d_{i-1}$.
\STATE Compute $FC(X)$ based on (\ref{first-column}), whose elements are $C_1,\cdots,C_\alpha$. For any $\forall~1\leq i\leq \alpha$, define $D_i=\mathfrak{N}_R(DCR(X,C_i))$.
\STATE Repeating step 1,2 on $D_i~(i=1,2,\cdots,\alpha)$, then we can get $ST(D_i)=\sum_{j=1}^{\xi_i}\lambda_{ij}E_{ij}$. \\
       By (\ref{sum-choice}), (\ref{column-sign}), (\ref{permutation}), for any $1\leq i\leq \alpha, 1\leq j\leq d$, compute
       $\Lambda_i, \epsilon_i, \pi_{ij}$.
\RETURN $\sum_{i=1}^\alpha~\sum_{j\in\Lambda_i}\textmd{sign}(\pi_{ij})^d(-1)^{\epsilon_i}\lambda_{ij}~~C_i\circ_R E_{ij}.$
\end{algorithmic}
\end{algorithm}

\bx Let $X=\begin{array}{ccc}
                \a_1 & \a_8 & \a_9 \\
                \a_2 & \a_5 & \a_7 \\
                \a_3 & \a_4 & \a_6 \\
              \end{array}$, and $d=n=3$. In algorithm \ref{a algorithm of compututation SB in the multilinear case},

\textbf{Step 1}, all vector variables in $X$ are $\{\a_1,\a_2,\a_3,\a_4,\a_5,\a_6,\a_7,\a_8,\a_9\}$, so by definition,
\[\mathcal{U}_1=\{\a_1\};~~~\mathcal{U}_2=\{\a_1,\a_2,\a_3,\a_4\};~~~\mathcal{U}_3=\{\a_1,\a_2,\a_3,\a_4,\a_5,\a_6,\a_7\}.\]

\textbf{Step 2}, by (\ref{first-column}), all elements in $FC(X)$ are:
\[C_1=\begin{array}{c}
    \a_1 \\
    \a_2 \\
    \a_3 \\
  \end{array},~~
C_2=\begin{array}{c}
    \a_1 \\
    \a_2 \\
    \a_4 \\
  \end{array},~~
C_3=\begin{array}{c}
    \a_1 \\
    \a_2 \\
    \a_6 \\
  \end{array},~~
C_4=\begin{array}{c}
    \a_1 \\
    \a_3 \\
    \a_5 \\
  \end{array},~~
C_5=\begin{array}{c}
    \a_1 \\
    \a_3 \\
    \a_7 \\
  \end{array},~~
C_6=\begin{array}{c}
    \a_1 \\
    \a_4 \\
    \a_5 \\
  \end{array},~~
C_7=\begin{array}{c}
    \a_1 \\
    \a_4 \\
    \a_7 \\
  \end{array}.\]
Then, we have
\[D_1=\begin{array}{cc}
                \a_4 & \a_6 \\
                \a_5 & \a_7 \\
                \a_8 & \a_9 \\
              \end{array},~~D_2=
              \begin{array}{cc}
                \a_3 & \a_6 \\
                \a_5 & \a_7 \\
                \a_8 & \a_9 \\
              \end{array},~~D_3=
              \begin{array}{cc}
                \a_3 & \a_4 \\
                \a_5 & \a_7 \\
                \a_8 & \a_9 \\
              \end{array},~~D_4=
              \begin{array}{cc}
                \a_2 & \a_7 \\
                \a_4 & \a_6 \\
                \a_8 & \a_9 \\
              \end{array},~~D_5=\begin{array}{cc}
                \a_2 & \a_5 \\
                \a_4 & \a_6 \\
                \a_8 & \a_9 \\
              \end{array},~~D_6=
              \begin{array}{cc}
                \a_2 & \a_7 \\
                \a_3 & \a_6 \\
                \a_8 & \a_9 \\
              \end{array},~~D_7=
              \begin{array}{cc}
                \a_2 & \a_5 \\
                \a_3 & \a_6 \\
                \a_8 & \a_9 \\
              \end{array}.\]

\textbf{Step 3}, compute $ST(D_i),~(i=1,2,\cdots,7)$. For $D_1$, all vector variables are $\{\a_4,$ $\a_5,$ $\a_6,$ $\a_7,$ $\a_8,$ $\a_9\}$.
$FC(D_1)=\{C_{11}=\begin{array}{c}
                \a_4 \\
                \a_5 \\
                \a_8 \\
              \end{array}\},$
so $D_{11}=\mathfrak{N}_R(DCR(D_i,C_{11}))=\begin{array}{c}
                \a_6 \\
                \a_7 \\
                \a_9 \\
              \end{array}.$
Thus $ST(D_{11})=D_{11}$, therefore $ST(D_1)=D_1$.

For the other six, the calculations are similar, so we omit the details, and only give the final results:
$ST(D_2)=D_2, ST(D_3)=D_3, ST(D_5)=D_5, ST(D_7)=D_7, ST(D_4)=
              \begin{array}{cc}
                \a_2 & \a_6 \\
                \a_4 & \a_7 \\
                \a_8 & \a_9 \\
              \end{array}-
              \begin{array}{cc}
                \a_2 & \a_4 \\
                \a_6 & \a_7 \\
                \a_8 & \a_9 \\
              \end{array},$ $ST(D_6)=
              \begin{array}{cc}
                \a_2 & \a_6 \\
                \a_3 & \a_7 \\
                \a_8 & \a_9 \\
              \end{array}-
              \begin{array}{cc}
                \a_2 & \a_3 \\
                \a_6 & \a_7 \\
                \a_8 & \a_9 \\
              \end{array}$.

\textbf{Step 4}, $\Lambda_1=\Lambda_2=\Lambda_3=\Lambda_5=\Lambda_7=\{1\},\Lambda_4=\{1,2\},\Lambda_6=\{2\}$,
$\epsilon_1=0,\epsilon_2=1,\epsilon_3=2,\epsilon_4=1,\epsilon_5=2,\epsilon_6=2,\epsilon_7=3$,
all $\pi_{ij}$ are identity permutation.

Finally, we have
\[ST(X)=\begin{array}{ccc}
                \a_1 & \a_4 & \a_6 \\
                \a_2 & \a_5 & \a_7 \\
                \a_3 & \a_8 & \a_9 \\
              \end{array}-
              \begin{array}{ccc}
                \a_1 & \a_3 & \a_6 \\
                \a_2 & \a_5 & \a_7 \\
                \a_4 & \a_8 & \a_9 \\
              \end{array}+
              \begin{array}{ccc}
                \a_1 & \a_3 & \a_4 \\
                \a_2 & \a_5 & \a_7 \\
                \a_6 & \a_8 & \a_9 \\
              \end{array}-
              \begin{array}{ccc}
                \a_1 & \a_2 & \a_7 \\
                \a_3 & \a_4 & \a_6 \\
                \a_5 & \a_8 & \a_9 \\
              \end{array}+
              \begin{array}{ccc}
                \a_1 & \a_2 & \a_4 \\
                \a_3 & \a_6 & \a_7 \\
                \a_5 & \a_8 & \a_9 \\
              \end{array}+\begin{array}{ccc}
                \a_1 & \a_2 & \a_5 \\
                \a_3 & \a_4 & \a_6 \\
                \a_7 & \a_8 & \a_9 \\
              \end{array}+
              \begin{array}{ccc}
                \a_1 & \a_2 & \a_6 \\
                \a_4 & \a_3 & \a_7 \\
                \a_5 & \a_8 & \a_9 \\
              \end{array}-
              \begin{array}{ccc}
                \a_1 & \a_2 & \a_3 \\
                \a_4 & \a_6 & \a_7 \\
                \a_5 & \a_8 & \a_9 \\
              \end{array}.\]
\ex

For the general case, the computation of $SB(X)$ are similar to the multi-linear case. Based on the analysis of multi-linear case, we only need to find the set $FC(X)$ in (\ref{first-column}). For the general case, the first column set of $SB(X)$ is denoted as $GFC(X)$.

Let $X$ be a Young tableau of degree $d$ and dimension $n$. Assume that the multiset of all the vector variables appeared in $X$ are
\be
\mathcal{M}(X)=\{\underbrace{\v_1,\cdots,\v_1}_{\alpha_1},\underbrace{\v_2,\cdots,\v_2}_{\alpha_2},\cdots,\underbrace{\v_k,\cdots,\v_k}_{\alpha_k}\}
:=\{\v_1^{\alpha_1},\v_2^{\alpha_2},\cdots,\v_k^{\alpha_k}\},
\ee
where $\alpha_i>0$ is the multiplicity of $\v_i$. Assume that $\v_1\prec\v_2\prec\cdots\prec\v_k$ lexicographically.

Let $\w\in\mathcal{M}(X)$, define
\be
InRow(X,\w):=\Big\{i~~\Big|~~1\leq i\leq d,~\w~\textmd{appears~in~the}~i\textmd{-th~row~of}~X\Big\}.
\ee

If $\begin{array}{c}
       \v_1^{\gamma_1} \\
       \v_2^{\gamma_2} \\
       \vdots \\
       \v_l^{\gamma_l} \\
     \end{array}\in GFC(X)$, then

(0). $l\leq k-(n-1)$. This is because that in the first column, the vector variable $\v_l$ with the highest lexicographical order satisfies $k-l\geq n-1$. This restriction can ensure that there are enough vector variables higher than $\v_l$ that can form a row in a term of $ST(X)$ with $\v_l$.

(1). $\gamma_1=\alpha_1$. This is because the lowest element can only appears in the left-upper of a Young tableau.

(2). $\max\{0,\alpha_2-\gamma_1\}\leq \gamma_2\leq \#(InRow(X_2,\v_2))$, where $X_2$ is the tableau obtained by deleting the rows containing $\v_1$ in $X$. In order to make sure that $\v_2$ and $\v_1$ do not lie in the same row in $X$, so $\v_2$ can only be chosen in $X_2$. Since in the choosing rows that containing $\v_1$, the number $\v_2$ can not surpass $\gamma_1$, so $\max\{0,\alpha_2-\gamma_1\}\leq \gamma_2$, which implies $\alpha_2-\gamma_2\leq \gamma_1$.

(3). $\max\{0,\alpha_3-\gamma_1-\gamma_2\} \leq \gamma_3 \leq \#(InRow(X_3,\v_3))$, where $X_3$ is the tableau obtained by deleting the rows containing $\v_2$ in $X_2$. In order to make sure that $\v_3$ and $\v_1,\v_2$ do not lie in the same row in $X$, so $\v_3$ can only be chosen in $X_3$. Since in the choosing rows that containing $\v_1,\v_2$, the number $\v_3$ can not surpass $\gamma_1+\gamma_2$, so $\max\{0,\alpha_3-\gamma_1-\gamma_2\}\leq \gamma_3$, which implies $\alpha_3-\gamma_3\leq \gamma_1+\gamma_2$.

$\cdots\cdots\cdots$

($i$). $\max\{0,\alpha_i-\sum_{j=1}^{i-1}\gamma_j\} \leq \gamma_i \leq \#(InRow(X_i,\v_i))$, where $X_i$ is the tableau obtained by deleting the rows containing $\v_{i-1}$ in $X_{i-1}$. In order to make sure that $\v_i$ and $\v_1,\v_2,\cdots,\v_{i-1}$ do not lie in the same row in $X$, so $\v_i$ can only be chosen in $X_i$. Since in the choosing rows that containing $\v_1,\v_2,\cdots,\v_{i-1}$, the number $\v_i$ can not surpass $\sum_{j=1}^{i-1}\gamma_j$, so $\max\{0,\alpha_i-\sum_{j=1}^{i-1}\gamma_j\}\leq \gamma_i$, which implies $\alpha_i-\gamma_i\leq \sum_{j=1}^{i-1}\gamma_j$.

$\cdots\cdots\cdots$

($l$). $\max\{0,\alpha_l-\sum_{j=1}^{l-1}\gamma_j\} \leq \gamma_l \leq \#(InRow(X_l,\v_l))$, where $X_l$ is the tableau obtained by deleting the rows containing $\v_{l-1}$ in $X_{l-1}$.

By the above analysis, we have

\bp
Let $X$ be a Young tableau of degree $d$ and dimension $n$. Then the first column set of $SB(X)$ is
\be \label{general-first-column}
GFC(X)=\Big\{\begin{array}{c}
       \v_1^{\gamma_1} \\
       \v_2^{\gamma_2} \\
       \vdots \\
       \v_l^{\gamma_l} \\
     \end{array}  ~~\Big|~~  l,\gamma_1,\gamma_2,\cdots,\gamma_l~\textmd{satisfy}~(0),(1),\cdots,(l)\Big\}.
\ee
\ep

The following procedure is similar to the multilinear case when $GFC(X)$ is achieved. Assume that $GFC(X)=\{C_1,C_2,\cdots,C_\alpha\}$. For any $1\leq i\leq \alpha$, denote $D_i=\mathfrak{N}_R(DCR(X,C_i))$, which is a Young tableau of dimension $(n-1)$. By induction, for any $1\leq i\leq \alpha$, we assume that $ST(D_i)$ has been computed, denoted as $ST(D_i)=\sum_{j=1}^{\xi_i}\lambda_{ij}E_{ij}$. For $1\leq i\leq \alpha$, define
\be \label{sum-choice1}
\Lambda_i=\Big\{1\leq j\leq \xi_i~~\Big|~~C_i\circ_R E_{ij}~\textmd{is~straight}\Big\}.
\ee
For any $1\leq i\leq \alpha, 1\leq j\leq d$, assume the $j$-th element of $C_i$ lies in the $q_{ij}$-th column of $X$, define
\be \label{column-sign1}
\epsilon_i=\sum_{j=1}^d (q_{ij}-1).
\ee
Denoting the first row of $X$ as $\t_1~\t_2~\cdots~\t_n$. For any $1\leq k\leq n$, let $\t_k$ appears in the $\omega_{ijk}$-th column of $C_i\circ_R E_{ij}$,
define the permutation $\pi_{ij}\in S_n$ as
\be \label{permutation1}
\pi_{ij}(k)=\omega_{ijk},~(\forall~k=1,2,\cdots,n).
\ee

In the above procedure, there may exist $1\leq i\neq j \leq \alpha$, such that $C_i=C_j$, while $D_i\neq D_j$, so the definition of (\ref{column-sign1}) and (\ref{permutation1}) will not make any confusion.

By the above analysis and the analysis of multilinear case, we have the following algorithm to compute $ST(X)$ in the general case.
\begin{algorithm}[H]
\caption{Compute $ST(X)$ in the general case}
\label{a algorithm of compututation SB in the general case}
\begin{algorithmic}[1]
\REQUIRE
A Young tableau $X$ of degree $d$ and dimension $n$
\ENSURE
$ST(X)$
\STATE Compute $\mathcal{M}(X)=\{\v_1^{\alpha_1},\v_2^{\alpha_2},\cdots,\v_k^{\alpha_k}\}$, where $\alpha_i>0$ and $\v_1\prec\v_2\prec\cdots\prec\v_k$ lexicographically.
\STATE Compute $GFC(X)$ based on (\ref{general-first-column}), whose elements are $C_1,\cdots,C_\alpha$. For any $\forall~1\leq i\leq \alpha$, define $D_i=\mathfrak{N}_R(DCR(X,C_i))$.
\STATE Repeating step 1,2 on $D_i~(i=1,2,\cdots,\alpha)$, then we can get $ST(D_i)=\sum_{j=1}^{\xi_i}\lambda_{ij}E_{ij}$. \\
       By (\ref{sum-choice1}), (\ref{column-sign1}), (\ref{permutation1}), for any $1\leq i\leq \alpha, 1\leq j\leq d$, compute
       $\Lambda_i, \epsilon_i, \pi_{ij}$.
\RETURN $\sum_{i=1}^\alpha~\sum_{j\in\Lambda_i}\textmd{sign}(\pi_{ij})^d(-1)^{\epsilon_i}\lambda_{ij}~~C_i\circ_R E_{ij}.$
\end{algorithmic}
\end{algorithm}

\bx Let $X=\begin{array}{ccc}
                     \a_1 & \a_2 & \a_6 \\
                     \a_2 & \a_3 & \a_4 \\
                     \a_2 & \a_4 & \a_5 \\
                   \end{array}$, then

\textbf{Step 1,} $\mathcal{M}(X)=\{\a_1,\a_2^3,\a_3, \a_4^2,\a_5,\a_6\}$, and $k=6,n=d=3$. So $l\leq 4$.

\textbf{Step 2,} by (\ref{general-first-column}),  $GFC(X)$ contains one element
$C_1=\begin{array}{c}
                     \a_1  \\
                     \a_2  \\
                     \a_2  \\
                   \end{array}.$
So we have
$D_1=\begin{array}{cc}
                     \a_2 & \a_6 \\
                     \a_3 & \a_4 \\
                     \a_4 & \a_5 \\
                   \end{array}.$

\textbf{Step 3,} computing $ST(D_1)$.
\[ST(D_1) = \begin{array}{cc}
                     \a_2 & \a_4 \\
                     \a_3 & \a_5 \\
                     \a_4 & \a_6 \\
                   \end{array}-\begin{array}{cc}
                     \a_2 & \a_4 \\
                     \a_3 & \a_4 \\
                     \a_5 & \a_6 \\
                   \end{array}-\begin{array}{cc}
                     \a_2 & \a_3 \\
                     \a_4 & \a_5 \\
                     \a_4 & \a_6 \\
                   \end{array}.\]

\textbf{Step 4,} $\Lambda_1=\{1,2,3\}$; $\varepsilon_1=0$ and $\pi_{11}=\pi_{12}=\pi_{13}=\textmd{id}$. Finally, we have
\[ST(X)=\begin{array}{ccc}
                     \a_1 & \a_2 & \a_4 \\
                     \a_2 & \a_3 & \a_5 \\
                     \a_2 & \a_4 & \a_6 \\
                   \end{array}-\begin{array}{ccc}
                     \a_1 & \a_2 & \a_4 \\
                     \a_2 & \a_3 & \a_4 \\
                     \a_2 & \a_5 & \a_6 \\
                   \end{array}-\begin{array}{ccc}
                     \a_1 & \a_2 & \a_3 \\
                     \a_2 & \a_4 & \a_5 \\
                     \a_2 & \a_4 & \a_6 \\
                   \end{array}.\]

\ex

\section{Testing of straightening algorithms}
\setcounter{equation}{0}

In this section, we will test the four straightening algorithms, i.e., classical straightening algorithm \verb"vw", White' straightening algorithm (we use algorithm C) \verb"white", Rota's straightening algorithm \verb"rota" and the straightening algorithm based on dual bracket \verb"db", by specific examples. Then we collect the data which can reflects the complexity of the straightening algorithms.

\verb"vw" and \verb"white" are the reduction procedures of Gr\"{o}bner bases. Given a non-straight bracket polynomial $F$, in these two straightening algorithms, one reduction by a Gr\"{o}bner base will be viewed as one step. Then we will get a new expression $F'$ of $F$ after each reduction, collect the number of terms of $F'$. Finally, we will get the data of total steps and the terms in each step of straightening $F$ by these two straightening algorithms.

In \verb"db", given a non-straight bracket polynomial $F=\sum_{i=1}^M \lambda_i T_i$, where each $T_i$ is a bracket monomial. In algorithm \ref{straightening algorithm based on dual bracket}, computing $q=SB(F)=\sum_{i=1}^M \lambda_i SB(T_i)$, which will be viewed as $M$ steps. And each computation of $SB(LT(q))$ viewed as one step. In each step, after combination, we will get a new expression $q'$ of $q$, collecting the number of terms of $q'$. Finally, we will get the data of total steps and the terms in each step of straightening $F$.

The texting data of \verb"vw", \verb"white" and \verb"db" will written as in the following form:
\be\ba{lll} \vspace{.1cm}
&& \textmd{[number~of~terms~in~the~first~step,~number~of~terms~in~the~second~step,~and~so~on],} \\
&& \textmd{~total~steps,~maximal~term,~total terms.}
\ea\ee
For instance, [14, 7, 2], 3, 14, 23 means: [14, 7, 2] is the collection of terms in each step, 3 is the total steps, 14 is the maximal terms and 23 is the total terms. However, when the degree are higher than 5, the number of terms in each step will be very large and the total step will be very large too, so we only collect some main data, such as ``total steps, maximal terms and total terms".

\verb"rota" is a little different from the above three straightening algorithms, it changes the straightening process into a solving of a system of linear equations, and the coefficient matrix is lower triangular with diagonal entries non-zero. So the number of non-zero entries in this coefficient matrix can coarsely reflect the complexity of Rota's straightening algorithm. Therefore, in \verb"rota", we collect the data of the number of non-zero entries in the coefficient matrix. So in the following, the collecting data about \verb"rota" keeps invariant when the number of rows and columns do not change.

In the following, we will present examples classified by the number of rows $r$ and columns $c$ of the bracket monomials or just say $r\times c$ bracket monomials. We only list a few number of examples when the case is quite simple, and list more when the case is representative.

(a). $3\times 3$ case: Since this case is quite simple, so we only collect the data of the following three bracket monomials.
\[(\textmd{i}):~\left[
     \begin{array}{ccc}
       \a_1 & \a_8 & \a_9 \\
       \a_2 & \a_6 & \a_7 \\
       \a_3 & \a_4 & \a_5 \\
     \end{array}
   \right],~~~(\textmd{ii}):~\left[
     \begin{array}{ccc}
       \a_1 & \a_6 & \a_9 \\
       \a_2 & \a_3 & \a_7 \\
       \a_4 & \a_5 & \a_8 \\
     \end{array}
   \right],~~~(\textmd{iii}):~\left[
     \begin{array}{ccc}
       \a_1 & \a_5 & \a_7 \\
       \a_2 & \a_6 & \a_8 \\
       \a_3 & \a_4 & \a_9 \\
     \end{array}
   \right].\]

\verb"vw":
\begin{eqnarray*}
(\textmd{i})  && [5, 3, 7, 11, 15, 19, 23, 27, 31, 32, 30, 32, 36, 40, 41, 45, 47, 49, 46, 50, 54, 53, 57, 56, 56, 53, 55, 54, 51, 51,   \\
     && 47, 44, 45, 41, 41, 41, 39, 38, 35, 34, 32, 30, 28, 25, 23, 19, 19, 18, 18, 15, 11, 11, 10, 10, 7, 3],56,56,1813; \\
(\textmd{ii}) && [5, 3, 7, 5, 5, 7, 7, 11, 9, 9, 10],11,11,78; \\
(\textmd{iii})&& [5, 9, 13, 17, 21, 25, 23, 27, 31, 32, 36, 34, 32, 30, 32, 32, 29, 28, 30, 29, 29, 26, 23, 23, 23],25,36,639.
\end{eqnarray*}

\verb"white":
\begin{eqnarray*}
(\textmd{i})    && [3, 5, 7, 9, 11, 12, 14, 17, 13, 9],10,17,100; \\
(\textmd{ii})   && [3, 5, 6, 9, 10],5,10,33; \\
(\textmd{iii})  && [4, 8, 11, 14, 17, 19, 22, 26, 30, 28, 32, 29, 29, 30, 30, 27, 27, 23, 24, 24, 25],21,32,479.
\end{eqnarray*}

\verb"rota": 96.

\verb"db":
\begin{eqnarray*}
(\textmd{i})   && [14, 7, 2],3,14,23; \\
(\textmd{ii})  && [12, 13, 9, 8, 10, 11, 5, 4, 3, 2],10,13,77; \\
(\textmd{iii}) && [13, 15, 18, 19, 19, 17, 18, 17, 16, 15, 18, 13, 13, 11, 10, 9, 8, 7, 7, 5, 4, 3, 2],23,19,277.
\end{eqnarray*}

(b). $4\times 3$ case: In this case, \verb"vw" will be eliminated from the text, while the other three straightening algorithm seems so easy, so we collect the data of the following three bracket monomials.
\[(\textmd{iv}):~\left[
     \begin{array}{ccc}
       \a_1 & \a_{11}      & \a_{12} \\
       \a_2 & \a_9\hfill   & \a_{10} \\
       \a_3 & \a_7\hfill   & \a_8\hfill \\
       \a_4 & \a_5\hfill   & \a_6\hfill \\
     \end{array}
   \right],~~~(\textmd{v}):~\left[
     \begin{array}{ccc}
       \a_1 & \a_8         & \a_{12} \\
       \a_2 & \a_7\hfill   & \a_{11} \\
       \a_3 & \a_6\hfill   & \a_{10} \\
       \a_4 & \a_5\hfill   & \a_9\hfill \\
     \end{array}
   \right],~~~(\textmd{vi}):~\left[
     \begin{array}{ccc}
       \a_1 & \a_{10}      & \a_{11} \\
       \a_2 & \a_7\hfill   & \a_{12} \\
       \a_3 & \a_6\hfill   & \a_8\hfill \\
       \a_4 & \a_5\hfill   & \a_9\hfill \\
     \end{array}
   \right].\]

\verb"vw": No results after 2 hours, the other methods return the answer in few minutes.

\verb"white":
\[\begin{tabular}{rllll} \vspace{.1cm}
  (\textmd{iv})~  & [3, 5, 7, 9, 11, \ldots, 87, 88, 89, 85, 81],  & 135, & 143, & 12459; \\ \vspace{.1cm}
  (\textmd{v})~   & [3, 5, 7, 9, 11, \ldots, 120, 121, 121, 121],  & 296, & 186, & 36849; \\
  (\textmd{vi})~  & [4, 7, 9, 11, 13, \ldots, 181, 182, 183, 181], & 258, & 220, & 42641. \\
\end{tabular}\]

\verb"rota": 2297.

\verb"db":
\[\begin{tabular}{rllll} \vspace{.1cm}
  (\textmd{iv})~  & [54, 65, 76, 78, \ldots, 4, 3, 2],        & 77,  & 76,  & 3266; \\ \vspace{.1cm}
  (\textmd{v})~   & [45, 54, 57, 49, \ldots, 6, 6, 4, 3, 2],  & 121, & 60,  & 4656; \\
  (\textmd{vi})~  & [55, 64, 73, 79, \ldots, 5, 4, 3, 2],     & 177, & 119, & 12564. \\
\end{tabular}\]

(c). $5\times 3$ case: In this case, \verb"white" and \verb"db" reflects more differences, so we collect the data of the following six bracket monomials.

\[\begin{tabular}{rrr} \vspace{.1cm}
  (\textmd{vii}):~$\left[
    \begin{array}{ccc}
      \a_1 & \a_{14}    & \a_{15} \\
      \a_2 & \a_{12}    & \a_{13} \\
      \a_3 & \a_{10}    & \a_{11} \\
      \a_4 & \hfill\a_8 & \hfill\a_9 \\
      \a_5 & \hfill\a_6 & \hfill\a_7 \\
    \end{array}
  \right]$,~~~ & (\textmd{viii}):~$\left[
    \begin{array}{ccc}
      \a_1 & \hfill\a_3  & \hfill\a_9 \\
      \a_2 & \hfill\a_7  & \hfill\a_8 \\
      \a_4 & \a_{10}     & \a_{15} \\
      \a_5 & \a_{11}     & \a_{13} \\
      \a_6 & \a_{12}     & \a_{14} \\
    \end{array}
  \right]$,~~~ & (\textmd{ix}):~$\left[
    \begin{array}{ccc}
      \a_1 & \hfill\a_6  & \a_{10} \\
      \a_2 & \hfill\a_7  & \a_{11} \\
      \a_3 & \hfill\a_8  & \a_{12} \\
      \a_4 & \hfill\a_9  & \a_{15} \\
      \a_5 & \a_{13}     & \a_{14} \\
    \end{array}
  \right]$, \\
  (\textmd{x}):~$\left[
    \begin{array}{ccc}
      \a_1 & \a_{10}     & \a_{15} \\
      \a_2 & \hfill\a_9  & \a_{14} \\
      \a_3 & \hfill\a_8  & \a_{13} \\
      \a_4 & \hfill\a_7  & \a_{12} \\
      \a_5 & \hfill\a_6  & \a_{11} \\
    \end{array}
  \right]$,~~~ & (\textmd{xi}):~$\left[
    \begin{array}{ccc}
      \a_1 & \hfill\a_6  & \a_{10} \\
      \a_2 & \hfill\a_7  & \a_{12} \\
      \a_3 & \hfill\a_8  & \a_{15} \\
      \a_4 & \hfill\a_9  & \a_{14} \\
      \a_5 & \a_{11}     & \a_{13} \\
    \end{array}
  \right]$,~~~ & (\textmd{xii}):~$\left[
    \begin{array}{ccc}
      \a_1 & \hfill\a_6  & \a_{13} \\
      \a_2 & \hfill\a_7  & \a_{15} \\
      \a_3 & \hfill\a_8  & \a_{12} \\
      \a_4 & \hfill\a_9  & \a_{14} \\
      \a_5 & \a_{10}     & \a_{11} \\
    \end{array}
  \right]$. \\
\end{tabular}\]

At this time, the data about the total steps, maximal term and total terms are:
\[\begin{tabular}{rrrrrrr} \vspace{.2cm}
                 & \verb"white":  &         &            & \verb"db":          &        &          \\
(\textmd{vii})   &~~~ 2052        &~~~1069  &~~~1324824  &\hspace{2cm} 267     &~~~237  &~~~29453  \\
(\textmd{viii})  &~~~ 571         &~~~422   &~~~171310   &\hspace{2cm} 334     &~~~195  &~~~35892  \\
(\textmd{ix})    &~~~ 155         &~~~143   &~~~14661    &\hspace{2cm} 112     &~~~61   &~~~4272   \\
(\textmd{x})     &~~~ 11020       &~~~3090  &~~~24345233 &\hspace{2cm} 1718    &~~~797  &~~~856398 \\
(\textmd{xi})    &~~~ 2161        &~~~846   &~~~1321347  &\hspace{2cm} 423     &~~~240  &~~~59161  \\
(\textmd{xii})   &~~~ 11581       &~~~4093  &~~~35537137 &\hspace{2cm} 1396    &~~~762  &~~~634952 \\
\end{tabular}\]

(d). $6\times 3$ case: In this case, the computation costs too much time, so we only collect the data of the following one bracket monomial.
\[(\textmd{xiii}):~\left[
    \begin{array}{ccc}
      \a_1 & \a_{17}    & \a_{18} \\
      \a_2 & \a_{15}    & \a_{16} \\
      \a_3 & \a_{13}    & \a_{14} \\
      \a_4 & \a_{11}    & \a_{12} \\
      \a_5 & \hfill\a_9 & \a_{10} \\
      \a_6 & \hfill\a_7 & \a_8 \\
    \end{array}
  \right].\]
\[\begin{tabular}{rl} \vspace{.1cm}
   \verb"white":~ & No results after 24 hours;  \\
   \verb"db":~    & [\textmd{total steps, maximal term, total terms}]~=~[1827595,~1624,~2037].    \\
\end{tabular}\]

(e). Turnbull-Young invariant polynomial \cite{ty}. This polynomial depicts the condition of ten points of three dimensional projective space lie in a common quadric surface.
Let $\sigma$ be a permutation about the set $\{2,3,4,5,6\}$, defined as
\[\sigma:~~~2\mapsto 3,~~3\mapsto 4,~~4\mapsto 5,~~5\mapsto 6,~~6\mapsto 2.\]
Let $S$ be the set of all permutations about $\{7,8,9\}$.
For any $0\leq j\leq 6$, define
\[\rho_j:~~~\a_j\mapsto \a_7,~~~\a_7 \mapsto \a_j,~~~\a_k\mapsto \a_k~~~(\forall~0\leq k\leq 9, k\neq j,7).\]
Setting
\begin{eqnarray*}
  H &=& \sum_{i=0}^4\sum_{\tau\in S}~[\a_0\a_1\a_{\sigma^i(2)}\a_{\sigma^i(3)}][\a_0\a_{\sigma^i(4)}\a_{\sigma^i(5)}\a_{\sigma^i(6)}]
        [\a_1\a_{\sigma^i(5)}\a_{\tau(7)}\a_{\tau(8)}] \\
    & & \hspace{3.15cm} [\a_{\sigma^i(2)}\a_{\sigma^i(6)}\a_{\tau(8)}\a_{\tau(9)}][\a_{\sigma^i(3)}\a_{\sigma^i(4)}\a_{\tau(7)}\a_{\tau(9)}], \\
  K &=& H-\sum_{j=0}^6\rho_i(H).
\end{eqnarray*}
The Turnbull-Young invariant polynomial is defined as $-\frac{1}{20}K$. By the above definition, the Turnbull-Young invariant polynomial contains 240 terms.

We using three straightening algorithms: \verb"vw", \verb"white" and~\verb"db" to straighten $K$, and collect the data. Our result shows that the straight expression of $K$ possesses 473 terms, the concrete expression in appendix I. While \cite{white1} get 138 terms about the straight expression of $K$.
The computation was on DELL workstation with Intel(R) Xeon(R) CPU E5-2630v3 @ 2.40 GHz RAM @ 32 G. The data are:
\begin{table}[H]
\centering
\caption{\bf data of Turnbull-Young invariant polynomial}
\begin{tabular}{ccccc}
  \hline \vspace{.1cm}
  ~~~     Method   ~~~&~~~ Total steps ~~~&~~~ Maxiaml term   ~~~&~~~ Total terms            ~~~&~~~ Time ~~~   \\ \vspace{.1cm}
  ~~~ \verb"vw"    ~~~&~~~ 20834       ~~~&~~~ 6457           ~~~&~~~ $1.0\times10^8$        ~~~&~~~ 11h  ~~~   \\ \vspace{.1cm}
  ~~~ \verb"white" ~~~&~~~ 9845        ~~~&~~~ 2841           ~~~&~~~ $1.9\times10^7$        ~~~&~~~ 6h   ~~~   \\ \vspace{.1cm}
  ~~~ \verb"db"    ~~~&~~~ 713         ~~~&~~~ 678            ~~~&~~~ $3.8\times10^5$        ~~~&~~~ 3h   ~~~   \\ \hline
\end{tabular}
\end{table}

(g). In the following, we tested more examples, and the data will be given in curves: the green one represents \verb"vw", the red one represents \verb"white", the black one represents \verb"rota" and the blue one represents \verb"db". The curves are collections of the date of total steps, total terms and the product of total steps and total terms.

In the multi-linear case, when the number of rows and columns are 3, and vector variables are $\a_1,\a_2,\cdots,\a_9$, there are totally 238 non-straight bracket monomials. We tested them all. While in the $4\times 3$ case, with vector variables $\a_1,\a_2,\cdots,\a_{12}$. The total number of non-straight bracket monomials is 1463, which is too hard to show in a graph, so we do not decided to test them all. Instead, we make some simple changes on the 238 $3\times 3$ non-straight bracket monomials (i.e., let
$\begin{array}{ccc}
       \f_1 & \f_2 & \f_3 \\
       \f_4 & \f_5 & \f_6 \\
       \f_7 & \f_8 & \f_9 \\
     \end{array}$ be a $3\times 3$ non-straight bracket monomial, then change it into $\begin{array}{ccc}
       \f_1 & \f_2 & \a_{12} \\
       \f_4 & \f_5 & \a_{11} \\
       \f_7 & \f_8 & \a_{10} \\
       \g_1 & \g_2 & \g_3    \\
     \end{array}$, where $\g_1,\g_2,\g_3$ is the sequence by sorting $\f_9,\f_6,\f_3$ into ascending form) and get 238 $4\times 3$ non-straight bracket monomials. We make a remark that the following graphs are just a collection of data of straightening algorithms, not depend on an fixed index. The resulting data (total terms, total steps), except \verb"rota" which is a uniform distribution, are fit for the normal distribution, with the mean $\mu$ and variance $\sigma$:

\verb"white":

\[\begin{tabular}{lll}
  \hline
               &~~~~ total terms                                  &~~~~ total steps \\
  $3\times 3$  &~~~~ $\sigma=84.6,\mu=127.3$                      &~~~~ $\sigma=6.5,\mu=5.9$ \\
  $4\times 3$  &~~~~ $\sigma=6.7\times 10^3,\mu=9.9\times 10^3$   &~~~~ $\sigma=90.1,\mu=87.1$ \\
  \hline
\end{tabular}\]

\verb"db":

\[\begin{tabular}{lll}
  \hline
               &~~~~ total terms                                  &~~~~ total steps \\
  $3\times 3$  &~~~~ $\sigma=81.6,\mu=80.8$                       &~~~~ $\sigma=10.0,\mu=6.4$ \\
  $4\times 3$  &~~~~ $\sigma=1.9\times 10^3,\mu=2.3\times 10^3$   &~~~~ $\sigma=50.6,\mu=37.6$ \\
  \hline
\end{tabular}\]

\verb"vw":

\[\begin{tabular}{lll}
  \hline
               &~~~~ total terms                                  &~~~~ total steps \\
  $3\times 3$  &~~~~ $\sigma=296.8,\mu=551.8$                     &~~~~ $\sigma=14.4,\mu=15.6$ \\
  \hline
\end{tabular}\]

(f). Complexity analysis in the multilinear case. Assume that the dimension of the bracket monomial is $n$, with degree $d$. Define
\be\ba{lll} \vspace{.1cm}
M_{n,d} &:=& \textmd{the~number~of~all~normal~Young~tableaux~of~degree}~d~\textmd{and~dimension}~n; \\ \vspace{.1cm}
N_{n,d} &:=& \textmd{the~number~of~all~straight~Young~tableaux~of~degree}~d~\textmd{and~dimension}~n; \\ \vspace{.1cm}
S_{n,d} &:=& \max\{\textmd{term~of}~SB(T)\mid T~\textmd{is~a~Young~tableau~of~degree}~d,~\textmd{dimension}~n\}; \\
b_n     &:=& C_{n+1}^{[(n+1)/2]},~\textmd{where}~$[(n+1)/2]$~\textmd{is~the~largest~integer~smaller~than}~$(n+1)/2$.
\ea\ee
Note that in the multilinear case, $M_{n,d}=\ds\frac{(dn)!}{d!(n!)^d}$ and $N_{n,d}=\ds\frac{(dn)!}{\prod_{i=1}^d\prod_{j=1}^n(i+j-1)}$ can be computed by the hook-length formula \cite{frame}. By Stirling's approximation \cite{well} $n!\approx \sqrt{2\pi n}\ds\left(\frac{n}{e}\right)^n$, we have
\be
M_{n,d} \approx d^{d(n-1)}n^{-\frac{1}{2}(d-1)}, \hspace{1cm} b_n \approx 2^n n^{-\frac{1}{2}},  \hspace{1cm}  S_{n,d}\leq (n!)^{d-1}\approx e^{-(d-1)(n-\frac{1}{2})}n^{(d-1)(n+\frac{1}{2})},
\ee
and
\be\ba{lll} \vspace{.1cm}
  N_{n,d} &=& \ds\frac{(dn)!2!3!\cdots(d-1)!}{n!(n+1)!\cdots(n+d-1)!} \\  \vspace{.1cm}
          &\approx& \ds\sqrt{\frac{d!n!}{(n+d-1)!}}\frac{d^{dn}n^{dn}2^23^3\cdots(d-1)^{d-1}}{n^n(n+1)^{n+1}\cdots(n+d-1)^{n+d-1}} \\
          &\approx& \ds\frac{2^23^3\cdots(d-1)^{d-1}d^dd^{dn-\frac{d}{2}+\frac{1}{4}}n^{dn-\frac{n}{2}+\frac{1}{4}}}{(n+1)^{n+1}\cdots
          (n+d-2)^{n+d-2}(n+d-1)^{n+d-1}(n+d-1)^{\frac{1}{2}(n+d-1)+\frac{1}{4}}}.
\ea\ee
Denote $h=\sum_{i=1}^n i\ln i$ and $S(i)=\sum_{j=1}^ij=\ds \frac{i(i+1)}{2}$, then
\begin{eqnarray*}
  h &=& \sum_{i=1}^n (S(i)-S(i-1))\ln i=-\sum_{i=1}^{n-1} S(i)(\ln (i+1)-\ln i) + S(n) \ln n \\
    &=& -\sum_{i=1}^{n-1} \frac{i(i+1)}{2} \left(\frac{1}{i}-\frac{1}{2i^2}+\frac{1}{3i^3}+O(\frac{1}{i^4})\right) + \frac{1}{2}n(n+1) \ln n \\
    &=& -\sum_{i=1}^{n-1} \left( \frac{i+1}{2}-\frac{1}{4}-\frac{1}{12i}   \right) + \frac{1}{2}n(n+1) \ln n +O(1) \\
    &=& -\frac{1}{4}n^2+ \frac{1}{2}(n^2+n+\frac{1}{6}) \ln n +O(1)
\end{eqnarray*}
The last equality is due to $\sum_{i=1}^{n-1}i^{-1}=\ln n+O(1)$. If $f=2^23^3\cdots(d-1)^{d-1}d^d$ and $g=(n+1)^{n+1}\cdots(n+d-2)^{n+d-2}(n+d-1)^{n+d-1}$, then
\begin{eqnarray*}
\ln f &=& -\frac{1}{4}d^2+ \frac{1}{2}(d^2+d+\frac{1}{6})\ln d +O(1), \\
\ln g &=& -\frac{1}{4}(n+d-1)^2+\frac{1}{4}n^2+\frac{1}{2}((n+d-1)^2+(n+d-1)+\frac{1}{6})\ln (n+d-1) \\
      & & -\frac{1}{2}(n^2+n+\frac{1}{6})\ln n  +O(1).
\end{eqnarray*}
Hence,
\be\ba{lll} \vspace{.1cm}
N_{n,d} &=      & \ds e^{\ln f-\ln g} d^{dn-\frac{d}{2}+\frac{1}{4}}n^{dn-\frac{n}{2}+\frac{1}{4}}(n+d-1)^{-\frac{1}{2}(n+d-1)-\frac{1}{4}} \\
        &\approx& \ds e^{\frac{1}{2}(dn-n-d)} d^{dn+\frac{1}{2}d^2+\frac{1}{3}} n^{dn+\frac{1}{2}n^2+\frac{1}{3}}(n+d-1)^{-\frac{1}{2}(n+d-1)(n+d+1)-\frac{1}{3}}.
\ea\ee

(f.1). \verb"vw". If $d=2$. Consider it in the worst case. Since in the first step, the input bracket monomial will returns $C_{n+1}^2$ terms. In the second step, we should choose the leaning term of the output of the first step under the row order, so we will get at most $C_{n+1}^2+C_{n+1}^3-1$ terms. After this, in the third step, we need to choose the leaning term of the output of the second step under the row order, so we will get at most $C_{n+1}^2+C_{n+1}^3+C_{n+1}^4-2$ terms, and so on. Finally, the total terms of \verb"vw" when $d$ is 2 is less than or equal to
\be\ba{lll}
 & & C_{n+1}^2+\left(C_{n+1}^2+C_{n+1}^3-1\right)+\cdots+\left(C_{n+1}^2+C_{n+1}^3+\cdots+C_{n+1}^{n}-(n-2)\right) \\
 &=& \ds 3(n+1)2^n-(n+1)^2-\frac{1}{2}(n-1)(n-2)=O(n2^n),
\ea\ee
and the total steps is less than or equal to $(n-1)$.

If $d\geq 3$, then in the first step, the worst case returns $b_n$ terms. Since $d\geq 3$, so we can not make sure the non-straight column will backward after each step. Therefore, in the second step, the total term is at most $2b_n-1$. Then in the third step, the total term is at most $3b_n-2$. Note that there are $M_{n,d}-N_{n,d}$ non-straight bracket monomials of degree $d$ and dimension $n$, so \verb"vw" needs at most $M_{n,d}-N_{n,d}$ steps. Finally, the total terms of \verb"vw" is less than or equal to
\be\ba{lll}
        & & b_n+2b_n-1+3b_n-2+\cdots+(M_{n,d}-N_{n,d})b_n-(M_{n,d}-N_{n,d}-1) \\
        &=& \ds\frac{1}{2}(M_{n,d}-N_{n,d})\big\{(M_{n,d}-N_{n,d}+1)(b_n-1)+2)\big\} \\
        &=& O(M_{n,d}^2b_n)=O(d^{2d(n-1)}n^{-(d-\frac{1}{2})} 2^n),
\ea\ee
and the total steps is less than or equal to $(M_{n,d}-N_{n,d})$.

(f.2). \verb"white". The worest case is similar to \verb"vw".

(f.3). \verb"rota". The number of non-zero elements of the coefficient matrix in Rota's straightening algorithm is less than or equal to
$\ds\frac{1}{2}N_{n,d}(N_{n,d}+1)$.
The algorithm of enumerating all straight bracket monomials of degree $d$ and dimension $n$ has complexity $N_{n,d}$. So the complexity of \verb"rota" is
\be
O(N_{n,d}^3)=O\left(\ds e^{\frac{3}{2}(dn-n-d)} d^{3dn+\frac{3}{2}d^2+1} n^{3dn+\frac{3}{2}n^2+1} (n+d-1)^{-\frac{3}{2}(n+d-1)(n+d+1)-1} \right).
\ee

(f.4). \verb"db". The first step returns at most $S_{n,d}$ terms. Then in the second step, we should choose the leading term of the result of first step under the negative column order, and compute its special bracket polynomial. So in the second step, we get at most $2S_{n,d}-1$ terms. Similarly, in the third step, we get at most $3S_{n,d}-2$ terms, and so on. The total step is at most $N_{n,d}$. Since at the final step of algorithm \ref{straightening algorithm based on dual bracket}, all the terms cancels. Therefore, the total terms of \verb"db" is less than or equal to
\be\ba{lll}
& & S_{n,d}+2S_{n,d}-1+3S_{n,d}-2+\cdots+(N_{n,d}-1)S_{n,d}-(N_{n,d}-2) \\
&=& \ds\frac{1}{2}\big((S_{n,d}-1)N_{n,d}^2-(S_{n,d}-3)N_{n,d}-1\big)=O(N_{n,d}^2S_{n,d}) \\
&=& O\left(\ds e^{-\frac{1}{2}d} d^{2dn+d^2+\frac{2}{3}} n^{3nd+n^2+\frac{1}{2}d-n+\frac{1}{6}} (n+d-1)^{-(n+d-1)(n+d+1)-\frac{2}{3}} \right).
\ea\ee
and the total steps is less than or equal to $N_{n,d}$.

In a conclusion, we get the following tables of the complexity estimation of different straightening algorithms.

When $n=3$, we have

\begin{table}[H]
\centering
\caption{\bf complexity estimation of different straightening algorithms when $n=3$}
\begin{tabular}{llll}
\hline \vspace{.15cm}
Methods                     ~~& Total~terms                                        ~~& Total~steps                          ~~& Maximal~steps                                      \\ \vspace{.15cm}
\verb"vw" and \verb"white"  ~~& $\ds e^{4d\ln d-d\ln 3}$                           ~~& $\ds e^{2d\ln d-\frac{1}{2}d\ln 3}$  ~~& $\ds e^{2d\ln d-\frac{1}{2}d\ln 3}$                \\ \vspace{.15cm}
\verb"rota"                 ~~& $\ds e^{3d+9d\ln 3-12\ln d}$                       ~~& $\ds e^{d+3d\ln 3-4\ln d}$           ~~&                                                    \\ \vspace{.1cm}
\verb"db"                   ~~& $\ds e^{-\frac{1}{2}d+\frac{19}{2}d\ln 3-8\ln d}$  ~~& $\ds e^{d+3d\ln 3-4\ln d}$           ~~& $\ds e^{-\frac{3}{2}d-4\ln d+\frac{13}{2}d\ln 3}$  \\ \hline
\end{tabular}
\end{table}

When $n=4$, we have

\begin{table}[H]
\centering
\caption{\bf complexity estimation of different straightening algorithms when $n=4$}
\begin{tabular}{llll}
\hline \vspace{.15cm}
Methods                     ~~& Total~terms                                        ~~& Total~steps                                       ~~& Maximal~steps                \\ \vspace{.15cm}
\verb"vw" and \verb"white"  ~~& $\ds e^{6d\ln d-2d\ln 2}$                          ~~& $\ds e^{3d\ln d-d\ln 2}$                          ~~& $\ds e^{3d\ln d-d\ln 2}$     \\ \vspace{.15cm}
\verb"rota"                 ~~& $\ds e^{\frac{9}{2}d+24d\ln 2-\frac{45}{2}\ln d}$  ~~& $\ds e^{\frac{3}{2}d+8d\ln 2-\frac{15}{2}\ln d}$  ~~&                              \\ \vspace{.1cm}
\verb"db"                   ~~& $\ds e^{-\frac{1}{2}d+\frac{31}{2}d\ln 3-15\ln d}$ ~~& $\ds e^{\frac{3}{2}d+8d\ln 2-\frac{15}{2}\ln d}$  ~~& $\ds e^{-2d-\frac{17}{2}\ln d+17d\ln 2}$  \\ \hline
\end{tabular}
\end{table}

In the general case, we have
\begin{table}[H]
\centering
\caption{\bf complexity estimation of total terms of different straightening algorithms}
\begin{tabular}{ccc}
    \hline \vspace{.15cm}
    ~~~~Methods     ~~~~~~&~~~~~~ Total~terms  \\ \vspace{.15cm}
    ~~~~\verb"vw" and \verb"white"  ~~~~~~&~~~~~~ $\left\{
                                    \begin{array}{ll} \vspace{.1cm}
                                      n2^n, & \hbox{~if~$d=2$;} \\
                                      d^{2d(n-1)}n^{-(d-\frac{1}{2})} 2^n, & \hbox{~if~$d>2$.}
                                    \end{array}
                                  \right.$~~~~~~   \\ \vspace{.15cm}
~~~~\verb"rota" ~~~~~~&~~~~~~ $\ds e^{dn-n-d} d^{2dn+d^2+\frac{2}{3}} n^{2dn+n^2+\frac{2}{3}} (n+d-1)^{-(n+d-1)(n+d+1)-\frac{2}{3}}  $~~~~~~   \\ \vspace{.1cm}
~~~~\verb"db"   ~~~~~~&~~~~~~ $\ds e^{-\frac{1}{2}d} d^{2dn+d^2+\frac{2}{3}} n^{3nd+n^2+\frac{1}{2}d-n+\frac{1}{6}} (n+d-1)^{-(n+d-1)(n+d+1)-\frac{2}{3}}$~~~~~~     \\ \hline
\end{tabular}
\end{table}

\begin{table}[H]
\centering
\caption{\bf complexity estimation of total steps of different straightening algorithms}
\begin{tabular}{ccc}
    \hline \vspace{.15cm}
    ~~~~Methods     ~~~~~~&~~~~~~ Total~steps  \\ \vspace{.15cm}
    ~~~~\verb"vw" and \verb"white"  ~~~~~~&~~~~~~ $\left\{
                                    \begin{array}{ll} \vspace{.1cm}
                                      n-1, & \hbox{~if~$d=2$;} \\
                                      d^{d(n-1)}n^{-\frac{1}{2}(d-1)}, & \hbox{~if~$d>2$.}
                                    \end{array}
                                  \right.$~~~~~~   \\ \vspace{.15cm}
~~~~\verb"rota" ~~~~~~&~~~~~~ $\ds e^{\frac{1}{2}(dn-n-d)} d^{dn+\frac{1}{2}d^2+\frac{1}{3}} n^{dn+\frac{1}{2}n^2+\frac{1}{3}}(n+d-1)^{-\frac{1}{2}(n+d-1)(n+d+1)-\frac{1}{3}}$ ~~~~~~ \\ \vspace{.1cm}
~~~~\verb"db"   ~~~~~~&~~~~~~ $\ds e^{\frac{1}{2}(dn-n-d)} d^{dn+\frac{1}{2}d^2+\frac{1}{3}} n^{dn+\frac{1}{2}n^2+\frac{1}{3}}(n+d-1)^{-\frac{1}{2}(n+d-1)(n+d+1)-\frac{1}{3}}$~~~~~~     \\ \hline
\end{tabular}
\end{table}

\begin{table}[H]
\centering
\caption{\bf complexity estimation of maximal term of different straightening algorithms}
\begin{tabular}{ccc}
    \hline \vspace{.15cm}
    ~~~~Methods     ~~~~~~&~~~~~~ Maximal~term  \\ \vspace{.15cm}
    ~~~~\verb"vw" and \verb"white"  ~~~~~~&~~~~~~ $\left\{
                                    \begin{array}{ll} \vspace{.1cm}
                                      2^n, & \hbox{~if~$d=2$;} \\
                                      d^{d(n-1)}n^{-\frac{1}{2}d} 2^n, & \hbox{~if~$d>2$.}
                                    \end{array}
                                  \right.$~~~~~~   \\ \vspace{.1cm}
~~~~\verb"db"   ~~~~~~&~~~~~~ $\ds e^{-\frac{1}{2}n(d-1)} d^{dn+\frac{1}{2}d^2+\frac{1}{3}} n^{2dn+\frac{1}{2}n^2+\frac{1}{2}d-n-\frac{1}{6} } (n+d-1)^{-\frac{1}{2}(n+d-1)(n+d+1)-\frac{1}{3}}$~~~~~~     \\ \hline
\end{tabular}
\end{table}

\section{Conclusion}
\setcounter{equation}{0}

In this paper, we present another new straightening algorithm which is fairly similar to Rota's straightening algorithm, while do not has its defect. The examples and the complexity analysis shown in section 5 tell us that this new straightening algorithm is much better than the others when the degree and dimension are higher.

\section{Appendix I: straightening expression of Turnbull-Young invariant polynomial}

In the following, for convenience, we use the subscript $i$ of $\a_i$ to represent it, where $i=0,1,\cdots,9$. The integers outside the brackets refer to coefficients. Then the straightening expression of Turnbull-Young invariant polynomial equals to the following bracket polynomial times $-\frac{1}{20}$.

\[\begin{tabular}{rrr}
 +2[0134][0246][1367][2589][5789] & +2[0124][0236][1567][3589][4789] & +2[0125][0126][3457][3489][6789] \\
 +2[0123][0247][1467][3589][5689] & +2[0123][0157][2467][3489][5689] & +2[0124][0237][1457][3689][5689] \\
 +2[0124][0235][1357][4689][6789] & +2[0124][0236][1357][4589][6789] & +2[0136][0247][1457][2589][3689] \\
 +2[0124][0345][1367][2589][6789] & +2[0134][0235][1467][2589][6789] & +2[0134][0237][1467][2589][5689] \\
 +2[0135][0347][1467][2589][2689] & +2[0135][0246][1267][3489][5789] & +2[0124][0257][1357][3689][4689] \\
 +2[0134][0234][1567][2689][5789] & +2[0134][0146][2357][2689][5789] & +2[0134][0247][1357][2689][5689] \\
 +2[0123][0245][1347][5689][6789] & +2[0125][0236][1457][3489][6789] & +2[0123][0256][1457][3689][4789] \\
 +2[0124][0234][1367][5689][5789] & +2[0134][0147][2367][2589][5689] & +2[0123][0167][2467][3589][4589] \\
 +2[0135][0235][1267][4689][4789] & +2[0136][0237][1257][4589][4689] & +2[0145][0257][1267][3489][3689] \\
 +2[0135][0257][1367][2489][4689] & +2[0123][0256][1467][3489][5789] & +2[0136][0257][1267][3489][4589] \\
 +2[0126][0347][1357][2589][4689] & +2[0125][0137][2347][4689][5689] & +2[0124][0127][3467][3589][5689] \\
 +2[0125][0345][1467][2689][3789] & +2[0126][0357][1467][2489][3589] & +2[0126][0246][1347][3589][5789] \\
 +2[0123][0124][3467][5689][5789] & +2[0125][0347][1347][2689][5689] & +2[0123][0246][1467][3589][5789] \\
 +2[0124][0126][3567][3589][4789] & +2[0124][0136][2356][4789][5789] & +2[0124][0137][2467][3589][5689] \\
 +2[0123][0137][2467][4589][5689] & +2[0125][0346][1357][2689][4789] & +2[0123][0156][2467][3589][4789] \\
 +2[0135][0267][1467][2489][3589] & +2[0126][0137][2457][3589][4689] & +2[0124][0347][1567][2589][3689] \\
 +2[0124][0237][1357][4689][5689] & +2[0134][0136][2467][2589][5789] & +2[0125][0345][1367][2689][4789] \\
 +2[0123][0235][1567][4689][4789] & +2[0125][0245][1367][3489][6789] & +2[0135][0236][1257][4689][4789] \\
 +2[0135][0357][1467][2489][2689] & +2[0125][0345][1367][2489][6789] & +2[0124][0145][2367][3689][5789] \\
 +2[0124][0134][2567][3589][6789] & +2[0134][0237][1257][4689][5689] & +2[0136][0246][1257][3589][4789] \\
 +2[0134][0135][2567][2689][4789] & +2[0123][0235][1467][4689][5789] & +2[0126][0236][1457][3589][4789] \\
 +2[0123][0256][1367][4589][4789] & +2[0124][0356][1357][2489][6789] & +2[0136][0346][1457][2589][2789] \\
 +2[0123][0147][2367][4589][5689] & +2[0124][0235][1567][3689][4789] & +2[0124][0157][2367][3489][5689] \\
 +2[0134][0145][2367][2589][6789] & +2[0123][0346][1567][2589][4789] & +2[0126][0136][2357][4589][4789] \\
 +2[0125][0146][2357][3489][6789] & +2[0123][0156][2467][3489][5789] & +2[0124][0246][1357][3689][5789] \\
 +2[0134][0137][2457][2689][5689] & +2[0126][0126][3457][3589][4789] & +2[0126][0137][2367][4589][4589] \\
 +2[0123][0146][2457][3689][5789] & +2[0135][0247][1367][2589][4689] & +2[0135][0235][1467][2689][4789] \\
 +2[0134][0247][1267][3589][5689] & +2[0135][0236][1457][2489][6789] & +2[0134][0245][1357][2689][6789] \\
 +2[0125][0136][2467][3589][4789] & +2[0123][0467][1467][2589][3589] & +2[0134][0156][2367][2589][4789] \\
 +2[0125][0345][1347][2689][6789] & +2[0136][0246][1357][2589][4789] & +2[0135][0246][1457][2689][3789] \\
 +2[0124][0234][1357][5689][6789] & +2[0124][0247][1367][3589][5689] & +2[0134][0157][2357][2689][4689] \\
 +2[0246][0135][1367][2489][5789] & +2[0135][0267][1267][3489][4589] & +2[0124][0347][1357][2689][5689] \\
 +2[0124][0135][2457][3689][6789] & +2[0124][0135][2347][5689][6789] & +2[0125][0135][2467][3489][6789] \\
 +2[0134][0157][2367][2489][5689] & +2[0125][0237][1467][3589][4689] & +2[0134][0257][1257][3689][4689] \\
 +2[0134][0236][1467][2589][5789] & +2[0124][0356][1357][2689][4789] & +2[0135][0246][1467][2589][3789] \\
 +2[0124][0237][1367][4589][5689] & +2[0125][0137][2467][3489][5689] & +2[0126][0127][3467][3589][4589] \\
 +2[0124][0125][3467][3689][5789] & +2[0135][0237][1457][2689][4689] & +2[0124][0235][1346][5789][6789] \\
 \end{tabular}\]
\[\begin{tabular}{rrr}
 +2[0124][0256][1457][3689][3789] & +2[0124][0367][1467][2589][3589] & +2[0145][0367][1367][2489][2589] \\
 +2[0135][0134][2457][2689][6789] & +2[0134][0167][2589][2367][4589] & +2[0125][0267][1467][3489][3589] \\
 +2[0134][0235][1456][2789][6789] & +2[0134][0256][1357][2689][4789] & +2[0124][0135][2345][6789][6789] \\
 +2[[0125]0246][1367][3589][4789] & +2[0134][0356][1467][2589][2789] & +2[0146][0247][1257][3589][3689] \\
 +2[0125][0236][1347][4689][5789] & +2[0123][0125][3457][4689][6789] & +2[0136][0247][1347][2589][5689] \\
 +2[0125][0247][1367][3489][5689] & +2[0134][0136][2457][2589][6789] & +2[0124][0135][2356][4789][6789] \\
 +2[0134][0246][1257][3689][5789] & +2[0123][0134][2456][5789][6789] & +2[0124][0235][1456][3789][6789] \\
 +2[0134][0235][1246][5789][6789] & +2[0135][0246][1267][3589][4789] & +2[0134][0357][1457][2689][2689] \\
 +2[0134][0345][1567][2689][2789] & +2[0126][0347][1367][2589][4589] & +2[0135][0237][1247][4689][5689] \\
 +2[0134][0257][1467][2589][3689] & +2[0124][0136][2457][3689][5789] & +2[0134][0236][1457][2689][5789] \\
 +2[0135][0246][1247][3689][5789] & +2[0135][0247][1257][4689][3689] & +2[0134][0367][1467][2589][2589] \\
 +2[0126][0247][1467][3589][3589] & +2[0124][0267][1367][3589][4589] & +2[0125][0347][1457][2689][3689] \\
 +2[0123][0267][1367][4589][4589] & +2[0126][0247][1357][3589][4689] & +2[0134][0236][1247][5689][5789] \\
 +2[0134][0234][1257][5689][6789] & +2[0124][0237][1567][3589][4689] & +2[0126][0357][1357][2489][4689] \\
 +2[0123][0257][1357][4689][4689]  & $-$2[0124][0235][1345][6789][6789] & $-$2[0124][0136][2346][5789][5789]  \\
$-$2[0125][0246][1357][3489][6789] & $-$2[0123][0246][1567][3589][4789] & $-$2[0124][0245][1567][3689][3789] \\
$-$2[0134][0234][1567][2589][6789] & $-$2[0135][0367][1467][2489][2589] & $-$2[0123][0137][2457][4689][5689] \\
$-$2[0134][0156][2357][2689][4789] & $-$2[0124][0136][2467][3589][5789] & $-$2[0134][0237][1267][4589][5689] \\
$-$2[0134][0247][1367][2589][5689] & $-$2[0135][0237][1467][2589][4689] & $-$2[0125][0147][2367][3489][5689] \\
$-$2[0134][0147][2357][2689][5689] & $-$2[0124][0346][1357][2589][6789] & $-$2[0134][0235][1567][2689][4789] \\
$-$2[0124][0245][1357][3689][6789] & $-$2[0123][0267][1467][3589][4589] & $-$2[0135][0246][1257][3689][4789] \\
$-$2[0126][0127][3457][3589][4689] & $-$2[0126][0347][1467][2589][3589] & $-$2[0134][0246][1357][2689][5789] \\
$-$2[0126][0346][1357][2589][4789] & $-$2[0134][0234][1267][5689][5789] & $-$2[0125][0347][1357][2689][4689] \\
$-$2[0124][0236][1457][3689][5789] & $-$2[0135][0245][1347][2689][6789] & $-$2[0134][0236][1257][4689][5789] \\
$-$2[0125][0235][1467][3489][6789] & $-$2[0134][0235][1256][4789][6789] & $-$2[0124][0234][1567][3589][6789] \\
$-$2[0123][0156][2367][4789][4589] & $-$2[0125][0367][1467][2489][3589] & $-$2[0123][0245][1567][3689][4789] \\
$-$2[0134][0256][1367][2589][4789] & $-$2[0123][0126][3567][4589][4789] & $-$2[0124][0345][1367][2689][5789] \\
$-$2[0123][0256][1357][4689][4789] & $-$2[0124][0135][2567][3689][4789] & $-$2[0124][0245][1367][3689][5789] \\
$-$2[0134][0356][1457][2689][2789] & $-$2[0135][0257][1467][2489][3689] & $-$2[0124][0256][1357][3689][4789] \\
$-$2[0124][0235][1356][4789][6789] & $-$2[0134][0235][1267][4589][6789] & $-$2[0135][0236][1267][4589][4789] \\
$-$2[0124][0236][1457][3589][6789] & $-$2[0135][0247][1457][2689][3689] & $-$2[0134][0136][2567][2589][4789] \\
$-$2[0134][0246][1357][2589][6789] & $-$2[0123][0145][2467][3689][5789] & $-$2[0124][0234][1567][3689][5789] \\
$-$2[0125][0346][1347][2689][5789] & $-$2[0126][0357][1367][2489][4589] & $-$2[0134][0267][1467][2589][3589] \\
$-$2[0124][0146][2357][3589][6789] & $-$2[0136][0236][1257][4589][4789] & $-$2[0125][0136][2367][4589][4789] \\
$-$2[0124][0134][2356][5789][6789] & $-$2[0126][0246][1357][3589][4789] & $-$2[0135][0257][1267][3489][4689] \\
$-$2[0135][0246][1357][2489][6789] & $-$2[0124][0237][1467][3589][5689] & $-$2[0125][0145][2367][3489][6789] \\
$-$2[0134][0256][1467][2589][3789] & $-$2[0134][0235][1247][5689][6789] & $-$2[0124][0356][1367][2589][4789] \\
$-$2[0124][0135][2467][3689][5789] & $-$2[0124][0347][1367][2589][5689] & $-$2[0134][0237][1247][5689][5689] \\
$-$2[0125][0125][3467][3489][6789] & $-$2[0146][0257][1257][3489][3689] & $-$2[0145][0247][1267][3589][3689] \\
$-$2[0124][0137][2457][3689][5689] & $-$2[0126][0237][1457][3589][4689] & $-$2[0134][0247][1257][3689][5689] \\
$-$2[0136][0347][1457][2589][2689] & $-$2[0123][0347][1567][2589][4689] & $-$2[0125][0346][1457][2689][3789] \\
$-$2[0134][0235][1467][2689][5789] & $-$2[0124][0357][1357][2689][4689] & $-$2[0125][0347][1367][2589][4689] \\
$-$2[0134][0136][2457][2689][5789] & $-$2[0123][0257][1457][3689][4689] & $-$2[0123][0146][2456][3789][5789] \\
$-$2[0136][0257][1257][3489][4689] & $-$2[0125][0346][1367][2589][4789] & $-$2[0134][0256][1257][3689][4789] \\
$-$2[0134][0256][1257][3489][6789] & $-$2[0123][0136][2467][4589][5789] & $-$2[0126][0247][1457][3589][3689] \\
$-$2[0126][0247][1367][3589][4589] & $-$2[0135][0247][1247][3689][5689] & $-$2[0134][0237][1457][2689][5689] \\
$-$2[0124][0247][1357][3689][5689] & $-$2[0126][0137][2467][3589][4589] & $-$2[0134][0145][2357][2689][6789] \\
$-$2[0123][0145][2346][5789][6789] & $-$2[0136][0247][1357][2589][4689] & $-$2[0126][0137][2357][4589][4689] \\
$-$2[0125][0126][3467][3589][4789] & $-$2[0145][0357][1367][2489][2689] & $-$2[0124][0136][2357][4589][6789] \\
$-$2[0125][0136][2457][3489][6789] & $-$2[0135][0267][1367][2489][4589] & $-$2[0126][0247][1347][3589][5689] \\
$-$2[0135][0236][1247][4689][5789] & $-$2[0124][0357][1467][2589][3689] & $-$2[0123][0245][1367][4689][5789] \\
$-$2[0136][0246][1347][2589][5789] & $-$2[0124][0136][2367][4589][5789] & $-$2[0125][0235][1347][4689][6789] \\
$-$2[0125][0237][1347][4689][5689] & $-$2[0125][0127][3467][3489][5689] & $-$2[0134][0245][1267][3589][6789] \\
 \end{tabular}\]
\[\begin{tabular}{rrr}
$-$2[0123][0345][1567][2689][4789] & $-$2[0124][0236][1357][4689][5789] & $-$2[0125][0246][1347][3589][6789] \\
$-$2[0123][0135][2457][4689][6789] & $-$2[0135][0236][1467][2489][5789] & $-$2[0135][0236][1457][2689][4789] \\
$-$2[0134][0146][2367][2589][5789] & $-$2[0134][0235][1245][6789][6789] & $-$2[0124][0235][1367][4689][5789] \\
$-$2[0126][0257][1467][3489][3589] & $-$2[0123][0256][1467][3589][4789] & $-$2[0124][0257][1457][3689][3689] \\
$-$2[0126][0136][2457][3589][4789] & $-$2[0135][0346][1467][2589][2789] & $-$2[0125][0136][2347][4689][5789] \\
$-$2[0125][0236][1467][3589][4789] & $-$2[0135][0245][1467][2689][3789] & $-$2[0123][0356][1467][2589][4789] \\
$-$2[0135][0245][1267][3689][4789] & $-$2[0136][0247][1257][3589][4689] & $-$2[0135][0247][1267][3489][5689] \\
$-$2[0123][0356][1457][2489][6789] & $-$2[0123][0167][2367][4589][4589] & $-$2[0134][0267][1367][2589][4589] \\
$-$2[0136][0246][1457][2589][3789] & $-$2[0124][0124][3567][3689][5789] & $-$2[0124][0346][1357][2689][5789] \\
$-$2[0134][0134][2567][2689][5789] & $-$2[0135][0246][1367][2589][4789] & $-$2[0124][0267][1467][3589][3589] \\
$-$2[0124][0246][1367][3589][5789] & $-$2[0134][0257][1357][2689][4689] & $-$2[0125][0237][1467][3489][5689] \\
$-$2[0124][0235][1347][5689][6789] & $-$2[0135][0237][1257][4689][4689] & $-$2[0134][0257][1367][2489][5689] \\
 +4[0126][0136][2347][4589][5789] & +4[0124][0136][2457][3589][6789] & +4[0123][0257][1467][3589][4689] \\
 +4[0123][0146][2356][4789][5789] & +4[0124][0137][2357][4689][5689] & +4[0123][0136][2567][4589][4789] \\
 +4[0123][0234][1456][5789][6789] & +4[0124][0246][1357][3589][6789] & +4[0123][0357][1467][2489][5689] \\
 +4[0123][0345][1457][2689][6789] & +4[0125][0135][2347][4689][6789] & +4[0134][0236][1257][4589][6789] \\
 +4[0136][0246][1247][3589][5789] & +4[0125][0346][1347][2589][6789] & +4[0134][0134][2567][2589][6789] \\
 +4[0134][0137][2567][2589][4689] & +4[0124][0357][1367][2589][4689] & +4[0123][0346][1456][2789][5789] \\
 +4[0134][0256][1267][3589][4789] & +4[0126][0347][1347][2589][5689] & +4[0124][0346][1367][2589][5789] \\
 +4[0136][0237][1247][4589][5689] & +4[0125][0347][1367][2489][5689] & +4[0123][0246][1456][3789][5789] \\
 +4[0123][0346][1457][2589][6789] & +4[0124][0235][1467][3589][6789] & +4[0123][0234][1467][5689][5789] \\
 +4[0123][0247][1457][3689][5689] & +4[0123][0345][1467][2689][5789] & +4[0123][0347][1457][2689][5689] \\
 +4[0123][0357][1467][2589][4689] & +4[0124][0134][2567][3689][5789] & +4[0135][0246][1347][2589][6789] \\
 +4[0123][0356][1457][2689][4789] & +4[0123][0245][1467][3689][5789] & +4[0123][0156][2357][4689][4789] \\
 +4[0124][0125][3456][3789][6789] & +4[0123][0456][1467][2589][3789] & +4[0124][0356][1367][2489][5789] \\
 +4[0123][0134][2567][4689][5789] & +4[0134][0267][1267][3589][4589] & +4[0123][0247][1347][5689][5689] \\
 +4[0123][0346][1467][2589][5789] & +4[0134][0257][1367][2589][4689] & +4[0123][0236][1456][4789][5789] \\
 +4[0123][0157][2457][3689][4689] & +4[0135][0237][1467][2489][5689] & +4[0124][0257][1467][3589][3689] \\
 +4[0134][0245][1267][3689][5789] & +4[0123][0147][2347][5689][5689]    & $-$4[0124][0135][2456][3789][6789] \\
$-$4[0134][0237][1567][2589][4689] & $-$4[0124][0357][1367][2489][5689] & $-$4[0123][0145][2367][4689][5789] \\
$-$4[0123][0156][2457][3689][4789] & $-$4[0135][0245][1267][3489][6789] & $-$4[0123][0146][2357][4689][5789] \\
$-$4[0123][0146][2567][3589][4789] & $-$4[0124][0135][2467][3589][6789] & $-$4[0136][0236][1247][4589][5789] \\
$-$4[0126][0137][2347][4589][5689] & $-$4[0134][0256][1267][3489][5789] & $-$4[0124][0367][1367][2589][4589] \\
$-$4[0136][0247][1247][3589][5689] & $-$4[0124][0137][2567][3589][4689] & $-$4[0123][0146][2346][5789][5789] \\
$-$4[0123][0257][1367][4589][4689] & $-$4[0124][0257][1367][3589][4689] & $-$4[0123][0356][1467][2489][5789] \\
$-$4[0134][0357][1467][2589][2689] & $-$4[0123][0136][2456][4789][5789] & $-$4[0126][0346][1347][2589][5789] \\
$-$4[0134][0246][1267][3589][5789] & $-$4[0123][0357][1457][2689][4689] & $-$4[0124][0247][1567][3589][3689] \\
$-$4[0123][0134][2457][5689][6789] & $-$4[0123][0135][2567][4689][4789] & $-$4[0135][0247][1367][2489][5689] \\
$-$4[0123][0236][1457][4689][5789] & $-$4[0123][0245][1467][3589][6789] & $-$4[0134][0137][2467][2589][5689] \\
$-$4[0134][0157][2367][2589][4689] & $-$4[0124][0125][3457][3689][6789] & $-$4[0123][0234][1457][5689][6789] \\
$-$4[0134][0257][1267][3589][4689] & $-$4[0125][0346][1357][2489][6789] & $-$4[0134][0236][1457][2589][6789] \\
$-$4[0125][0346][1367][2489][5789] & $-$4[0123][0157][2357][4689][4689] & $-$4[0123][0145][2456][3789][6789] \\
$-$4[0134][0135][2467][2589][6789] & +6[0123][0145][2345][6789][6789]   & +6[0123][0157][2367][4589][4689]  \\
 +6[0134][0257][1267][3489][5689] & +6[0123][0135][2467][4589][6789] & +6[0136][0247][1257][3489][5689] \\
 +6[0135][0236][1247][4589][6789] & +6[0123][0457][1457][2689][3689] & +6[0123][0145][2457][3689][6789] \\
 +6[0124][0135][2367][4689][5789] & +6[0123][0145][2356][4789][6789] & +6[0123][0146][2367][4589][5789] \\
 +6[0123][0246][1457][3589][6789] & +6[0123][0136][2457][4689][5789] & +6[0123][0145][2567][3689][4789] \\
 +6[0123][0126][3457][4589][6789] & +6[0123][0126][3456][4789][5789] & +6[0126][0346][1357][2489][5789] \\
 +6[0123][0247][1367][4589][5689]  & $-$6[0124][0137][2347][5689][5689] & $-$6[0125][0136][2347][4589][6789] \\
$-$6[0123][0257][1467][3489][5689] & $-$6[0123][0137][2567][4589][4689] & $-$6[0123][0157][2467][3589][4689] \\
$-$6[0123][0457][1467][2589][3689] & $-$6[0124][0126][3457][3589][6789] & $-$6[0123][0247][1357][4689][5689] \\
$-$6[0136][0246][1257][3489][5789] & $-$6[0124][0136][2357][4689][5789] & $-$6[0124][0124][3567][3589][6789] \\
$-$6[0126][0347][1357][2489][5689] & $-$6[0123][0456][1457][2689][3789] & $-$6[0135][0246][1247][3589][6789] \\
 \end{tabular}\]
\[\begin{tabular}{rrr}
$-$6[0123][0145][2367][4589][6789] & $-$6[0123][0147][2467][3589][5689] & $-$6[0123][0145][2467][3589][6789] \\
$-$6[0123][0347][1467][2589][5689] & $-$6[0124][0135][2357][4689][6789] & $-$6[0124][0134][2367][5689][5789] \\
$-$6[0123][0124][3457][5689][6789] & +8[0124][0125][3467][3589][6789]   & +8[0123][0234][1567][4589][6789] \\
 +8[0123][0125][3467][4689][5789]  & +8[0123][0127][3457][4689][5689]   & +8[0123][0246][1346][5789][5789] \\
 +8[0124][0136][2347][5689][5789]  & +8[0123][0245][1367][4589][6789]   & +8[0123][0156][2457][3489][6789] \\
 +8[0123][0126][3467][4589][5789]  & +8[0123][0145][2347][5689][6789]   & +8[0123][0146][2357][4589][6789] \\
 +8[0135][0246][1257][3489][6789]  & $-$8[0123][0136][2457][4589][6789] & $-$8[0123][0246][1356][4789][5789] \\
$-$8[0123][0246][1457][3689][5789] & $-$8[0123][0245][1346][5789][6789] & $-$8[0123][0234][1567][4689][5789] \\
$-$8[0123][0345][1467][2589][6789] & $-$8[0123][0346][1457][2689][5789] & $-$8[0123][0256][1457][3489][6789] \\
$-$8[0123][0235][1467][4589][6789] & +10[0123][0236][1467][4589][5789]  & +10[0123][0147][2567][3589][4689]  \\
+10[0124][0135][2367][4589][6789]   & $-$10[0123][0125][3467][4589][6789]  & $-$10[0123][0237][1467][4589][5689]\\
$-$10[0123][0246][1357][4589][6789] & $-$10[0123][0246][1367][4589][5789]  & $-$10[0123][0124][3567][4689][5789] \\
 +12[0123][0123][4567][4689][5789]  & $-$12[0123][0134][2567][4589][6789]  & $-$12[0123][0246][1347][5689][5789]  \\
$-$12[0123][0123][4567][4589][6789] & +14[0123][0246][1357][4689][5789]    & $-$14[0123][0126][3457][4689][5789] \\
+18[0123][0124][3567][4589][6789]   & +22[0123][0124][3456][5789][6789]    &  \\
 \end{tabular}\]

\section{Appendix II: straightening coefficients}
\setcounter{equation}{0}

In this section, we will consider the formula of straightening coefficients in bracket algebra. This have been studied in \cite{clausen2}, \cite{desarmenien}, \cite{huang}. Their proofs are based on Rota's straightening algorithm \cite{rota-straightening-formula} and a technique raised by Clausen \cite{clausen2}. Now we will give a proof about the formulas based on dual bracket.

\bn \label{coeff}
Let $S,T$ be two Young tableaux of degree $d$ and dimension $n$, define
\be
\mathcal{E}(S,T):=\ds\sum_{ \substack{\exists~\sigma_2,\cdots,\sigma_d\in S_n,~s.t., \\ T=\mathfrak{N}_C(S_{\emph{id},\sigma_2,\cdots,\sigma_d})} } \emph{sign}(\sigma_2)\cdots\emph{sign}(\sigma_d).
\ee
\en

This notation $T=\mathfrak{N}_C(S_{\emph{id},\sigma_2,\cdots,\sigma_d})$ is similar to the notation $S \leftrightarrow X\updownarrow T$ in \cite{huang} or $\leftidx{^C}(\leftidx{^R}S)$ in \cite{clausen2}, \cite{desarmenien}.
By definition \ref{db:def} and definition \ref{coeff}, we have

\bl \label{lemma1}
Let $Y$ be a Young tableau, assume that $DT(Y)=\sum \lambda_i Y_i$, then $\lambda_i=\mathcal{E}(Y,Y_i)$.
\el

As a corollary of lemma \ref{lemma1} and algorithm \ref{straightening algorithm based on dual bracket}, we have the following result about
straightening coefficients.

\bc
Let $Y$ be a non-straight bracket monomial and $Y=\sum \lambda_i S_i$ be the straight expression of $Y$ in the descending form under the negative column order, then
\be\ba{lll} \vspace{.1cm} \label{straightening coefficient}
\lambda_1=1, \\
\lambda_i=\mathcal{E}(Y,S_i)-\ds\sum_{j<i} \lambda_j\mathcal{E}(S_j,S_i),~~i\geq 2.
\ea\ee
\ec

From (\ref{straightening coefficient}), we find that the straightening coefficient $\lambda_i$ depends not only on $Y$, but also on $S_j,j<i$.

\end{document}